\documentclass[12pt]{amsart}
\usepackage[T1]{fontenc}
\usepackage{lmodern,amsmath,amssymb,amsthm,mathtools,microtype}

\usepackage[margin=4cm]{geometry}

\usepackage{xypic}
\usepackage{graphicx}
\usepackage{tikz}
\usetikzlibrary{arrows.meta}

\usepackage{latexsym}
\usepackage{epsfig}
\usepackage{amsfonts}
\usepackage{enumerate}
\usepackage{times,float}
\usepackage{graphics}

\usepackage[nodisplayskipstretch]{setspace}
\usepackage{enumitem,needspace}
\usepackage[dvipsnames]{xcolor}
\usepackage[pagebackref=true,colorlinks=true,linkcolor=MidnightBlue,citecolor=MidnightBlue,urlcolor=MidnightBlue]{hyperref}

\newtheorem{theorem}{Theorem}[section]
\newtheorem{proposition}[theorem]{Proposition}
\newtheorem{lemma}[theorem]{Lemma}
\newtheorem{corollary}[theorem]{Corollary}
\newtheorem*{externaltheorem}{Polynomial realization}
\theoremstyle{definition}
\theoremstyle{remark}\newtheorem{remark}[theorem]{Remark}
\numberwithin{equation}{section}
\newcommand{\PP}{\mathbb P}

\newcommand{\Q}{\mathbf Q}
\newcommand{\OO}{\mathcal O}

\newcommand{\cX}{\mathcal X}

\newcommand{\cZ}{\mathcal Z}
\DeclareMathOperator{\Tot}{Tot}

\title[Campana rational connectedness]{klt weak log Fano orbifolds are Campana rationally connected}

\author{Jiun-Cheng Chen}
\address{Department of Mathematics\\ Third General Building\\ National Tsing-Hua University\\ No. 101 Sec. 2 Kuang Fu Road\\ Hsinchu, 30043\\ Taiwan}
\email{jcchen@math.nthu.edu.tw}

\author[Hsian-Hua Tseng]{Hsian-Hua Tseng}
\address{Department of Mathematics\\ Ohio State University\\ 100 Math Tower, 231 West 18th Ave.\\ Columbus, OH 43210\\ USA}
\email{hhtseng@math.ohio-state.edu}
\date{\today}

\begin{document}

\begin{abstract}
Let \(X\) be a normal projective complex variety with a boundary divisor 
\(\Delta=\sum_{i=1}^s(1-1/m_i)D_i\) 
and integers $m_i$. We prove that, if $(X,\Delta)$ is klt and \(-(K_X+\Delta)\) is nef and big, then \((X,\Delta)\) is Campana rationally connected. 

When $(X,\Delta)$ is log smooth, we prove the theorem  by constructing a free Campana rational curve whose class lies in the interior of the nef cone of curve classes and whose contacts with \(D_i\) have order \(m_i\) and avoid deeper strata. We construct such a curve from a free Campana curve by attaching to it free rational curves which intersect positively with divisors that span the cone of effective divisors, and smooth the resulting comb. 

When $X$ is singular, we construct suitable rational curves on a log resolution of $(X,\Delta)$ and descend them to $X$ by appropriate  modifications and a relative minimal model program.
\end{abstract}
\maketitle
\tableofcontents

\setcounter{section}{-1}

\section{Introduction}\label{sec:intro}
We work over $\mathbb{C}$. Let \(X\) be a normal projective variety of positive dimension. Let 
\(D=\sum_{i=1}^sD_i\) be a divisor with distinct irreducible components.

\subsection{Result}\label{sec:result}

 Fix integers \(m_i\ge2\) and put
\[
\Delta=\sum_{i=1}^s\left(1-\frac{1}{m_i}\right)D_i.
\]
A \emph{Campana rational curve} on $(X,\Delta)$ is a nonconstant morphism \(f:\PP^1\to X\) whose image meets \(X\setminus D\) and such that
\begin{equation}\label{eq:campana}
\text{for }p\in\PP^1,\ 1\le i\le s, \text{ if } c_{p,i}:=\operatorname{ord}_p(f^*D_i)>0,
\text{ then } c_{p,i}\ge m_i.
\end{equation}
If \(f(p)\) belongs to several components of $D$, we impose this condition for each of them. We mark the points in \(f^{-1}(D)\). The {\em contact data} of $f: \PP^1\to X$ with $D$ consist of the vectors 
$$(c_{p,i})_{i=1}^s$$ 
at these markings. We say that \((X,\Delta)\) is \emph{Campana rationally connected} if there is a family of rational curves with fixed contact data for which the evaluation map 
$$(f,p,q)\mapsto(f(p),f(q))\in X\times X,$$ 
with \(p,q\) varying on the domain curve, is dominant. In particular, through two general points of \(X\) there is a rational curve in this family.

The main result of this paper is the following.
\begin{theorem}\label{thm:main_gen}
Suppose that $(X,\Delta)$ is a projective  klt pair with standard coefficients and \(-(K_X+\Delta)\) is nef and big. Then there exists a free Campana rational curve 
\(f:\PP^1\to X\) such that 
\begin{itemize}
\item 
$f_*[\PP^1]$ lies in the interior of \(\operatorname{Nef}_1(X)\);
\item 
$f$ avoids \(D_i\cap D_j\) for all \(i\ne j\); 
\item 
$\operatorname{ord}_p(f^*D_i)=m_i$ for all $p\in f^{-1}(D_i)$.
\end{itemize}
Consequently, \((X,\Delta)\) is Campana rationally connected.
\end{theorem}

Here \(\operatorname{Nef}_1(X)\subset N_1(X)_{\mathbb R}\) is the cone dual to the cone of  pseudoeffective divisors. By freeness of a Campana curve we mean the notion in \cite[Definition~4.2]{CLT} defined using its logarithmic normal complex with all boundary divisor contacts marked. The exact orders in the theorem refer to the specified free curve. The family of two-pointed curves satisfies the inequalities \eqref{eq:campana}.

Theorem \ref{thm:main_gen} proves \cite[Conjecture~1.2]{LT}, see also \cite[Conjecture 11.5(2)]{Cam}. For smooth log Fano orbifolds of Picard number $1$, see \cite[Corollary 1.8]{CT}. 

In the ample case, the rational connectedness assertion is Campana's conjecture as stated in \cite[Conjecture~5.8]{CLT}. Free Campana curves of possibly higher genera are constructed in \cite{LT}. Moreover, it is prove in \cite[Theorems~1.3 and~1.9]{LT} that, for a klt log Fano pair, rational connectedness follows from the existence of a very free rational curve in its smooth locus. 

Campana also conjectured \cite[Conjecture~11.5(1)]{Cam} that a smooth projective orbifold pair 
$(X, \Delta)$ over $\mathbb{C}$ is orbifold uniruled if and only if $K_X+\Delta$ is not pseudo-effective. To the best of our knowledge, this conjecture is still open. 

We divide our proof of Theorem \ref{thm:main_gen} into two parts. In Part \ref{part:smooth}, we establish the smooth case. In Part \ref{part:klt}, we use results proved in Part \ref{part:smooth} and additional arguments to establish the klt case.

As we were finishing writing up Part \ref{part:klt}, we noticed on the arXiv the paper \cite{Xu}, which contains related results in the case $D=\emptyset$.

\subsection{Acknowledgment}
In the process of developing the results presented here and writing this article, we benefited from the assistance of GPT-6 and 6.2 Astra. J.-C. C. is supported by the National Science and Technology Council (NSTC) of Taiwan under Grant No. NSTC-114-2115-M-007-004- and NSTC-115-2115-M-007-013-.

\part{Smooth case}\label{part:smooth}
In this part, we prove the log smooth case of Theorem \ref{thm:main_gen}.
\section{Overview: smooth case}
The main result of this part is
\begin{theorem}\label{thm:main}
Suppose that $(X,D)$ is log smooth and \(-(K_X+\Delta)\) is nef and big. Then there exists a free Campana rational curve 
\(f:\PP^1\to X\) such that 
\begin{itemize}
\item 
$f_*[\PP^1]$ lies in the interior of \(\operatorname{Nef}_1(X)\);
\item 
$f$ avoids \(D_i\cap D_j\) for all \(i\ne j\); 
\item 
$\operatorname{ord}_p(f^*D_i)=m_i$ for all $p\in f^{-1}(D_i)$.
\end{itemize}
\end{theorem}

\subsection{Proof strategy}
By \cite[Theorem~1]{Zhang}, the variety \(X\) is rationally connected. For such an \(X\), it is shown in \cite{CLT} that one free Campana rational curve whose class is interior to 
\(\operatorname{Nef}_1(X)\) suffices to obtain a Campana rational curve through two general points, see \cite[Definition~5.9 and Corollary~6.7]{CLT}. We first construct covering families of rational curves satisfying the required contact conditions\footnote{This also follows from the results of \cite{CT}.}. To construct curves of classes in the interior of \(\operatorname{Nef}_1(X)\), we carry out another construction of rational curves that provides some control of intersection with an additional choice of effective divisor.

Let $$\pi:\cX\to X$$ 
be the stack obtained by taking the \(m_i\)-th root along each \(D_i\). The identity 
$K_{\cX}\sim_{\mathbb Q}\pi^*(K_X+\Delta)$ lets us compute the $-(K_X+\Delta)$-degree as $-K_\cX$-degree on \(\cX\). Also, root construction implies that, for a nonconstant map $g:\PP^1\to \cX$ from an ordinary \(\PP^1\) to \(\cX\), if its image meets the complement of the root divisors, then the coarse\footnote{I.e. the contact orders of the map $\pi\mathbin{\!\scriptstyle\circ\!} f:\PP^1\to X$.} contact orders along \(D_i\) are divisible by \(m_i\). The construction of such maps starts with a complete curve $C$ of negative canonical degree on a root stack. For a smooth projective target, the {\em polynomial realization theorem} \cite[Theorem~4.2]{GX} gives a morphism defined near the zero section of a positive-degree line bundle on a finite cover $C'\to C$ of $C$. The morphism extends the given map on the zero section and is nonconstant along a general fiber. Every closed curve in the total space that dominates the base $C'$ meets the zero section. It follows that if the original complete curve $C$ lies in a given open subset of the target, a general affine fiber maps entirely into that subset. This is worked out in details in Lemma~\ref{lem:open-line}.

To apply this result to a root stack, we lift the complete curve, after finite base change, to a scheme of frames for the root line bundles. The frames are taken relative to line bundles on an abelian variety. The abelian variety and the torus of frames contribute $0$ to the degree against $c_1$ of the tangent bundle, and neither admits a nonconstant map from \(\mathbb A^1\). Consequently, the affine line just constructed has nonconstant projection to the coarse space. It already has the required divisibilities at its finite points. A power map ramified at infinity gives the remaining divisibilities when we complete the domain curve to \(\PP^1\), see Proposition~\ref{prop:incidence}. Applying this to lifts of general complete-intersection curves, we obtain a covering family. Variation of the member through a general point gives freeness on the stack, see Lemma~\ref{lem:covering}.

To control the numerical class of the curve, we need to consider its relation with a nonzero effective divisor \(E\) on \(X\). A covering curve has nonnegative \(E\)-degree, but its degree can be zero. Our next step is to exclude this possibility. We sketch the argument, see Proposition~\ref{prop:testing} for details. We first construct an auxiliary klt pair $(Y,B_Y)$, with 
$$q:Y\to X$$ 
a projective bundle. For a suitable integer \(m\ge3\), the fibers have dimension \(m-2\). On each fiber, we take \(m\) hyperplanes of weight \(1-1/m\). Their total degree \(m-1\) cancels the canonical degree of a line in the fiber. One of the hyperplanes varies with a section on \(X\). Choosing its line bundle in terms of \(E\), we obtain a boundary divisor \(B_Y\) satisfying
\[
K_Y+B_Y\sim_{\mathbb Q}-\varepsilon q^*E,\quad \varepsilon>0.
\]
Over the zero divisor of the section of $q$, the varying hyperplane passes through the common point of the affine coordinate hyperplanes. Under the nef and big hypothesis, the linear system of the section can have fixed components. We choose the section by a klt perturbation on \(X\) and check that the resulting pair on \(Y\) is klt, see Lemmas~\ref{lem:perturbation} and~\ref{lem:arrangement} for details.

We then  work with a log resolution $$\mu:\widetilde Y\to Y$$ 
and choose root orders so that the canonical class $K_{\mathcal{S}}$ of the resulting root stack $p:\mathcal{S}\to\widetilde{Y}$ is pulled back from \(-\varepsilon\mu^*q^*E+R\), where \(R\) is effective and $\mu$-exceptional, see (\ref{eq:exceptional}). Since \(R\) is exceptional, its intersection with a complete-intersection curve class pulled back from \(Y\) vanishes. A small ample perturbation therefore yields the negative degree needed to apply Lemma~\ref{lem:covering}. A general free rational curve upstairs has anticanonical degree at least $2$ and nonnegative degree against \(R\). The sum of these degrees is \(\varepsilon\) times the \(E\)-degree of the projected curve. This proves the strict positivity in Proposition~\ref{prop:testing}.

Finally, \(X\) is of Fano type, so the cone of its pseudoeffective divisors has finitely many effective generators. For each generator, we construct a free rational curve, namely a covering family of rational curves whose class intersects positively with this generator (and nonnegatively with all other generators). Pick a point \(x\in X\setminus D\) and form a map to $\cX$ from a comb by attaching suitable free curves from these families (``teeth'') to a $\PP^1$ (``handle'') mapping to $x$, see Figure \ref{fig:comb-smoothing} in Section \ref{sec:completion}. We then deform this map to obtain a map \(\PP^1\to \cX\), whose class is the sum of the classes of the free curves and is therefore in the interior, see Figure \ref{fig:campana-comb-smoothing} in Section \ref{sec:completion}. A general such deformation gives a map that avoids the deeper strata (i.e. intersections of two or more divisor components) and meets each root divisor transversely, so the contacts on \(X\) have exactly the orders \(m_i\). At each contact, the normal derivative identifies the tangent line of the domain curve with the normal line of the root divisor. The logarithmic normal complex downstairs is therefore quasi-isomorphic to the ordinary normal complex upstairs. We use this comparison in Section~\ref{sec:completion} to verify the freeness required by the two-point criterion.

We present two examples in Section \ref{sec:example}.

\section{An affine line in a smooth open target}\label{sec:open}

We will need an affine line through a point of a complete curve in a quasi-projective variety. We start with the following consequence of results from \cite{GX}.

\begin{externaltheorem}[Guo--Xie, {\cite[Proposition~3.2 and Theorem~4.2]{GX}}]
Let \(Y\) be a smooth projective variety and \(b:C\to Y\) a morphism from a smooth connected projective curve with \(d=-\deg b^*K_Y>0\). Then there exist
\begin{enumerate}
\item 
a finite surjective morphism \(\rho:C'\to C\) from a smooth connected projective curve,
\item 
a line bundle \(M\) on \(C'\) with \(\deg M>0\),
\item 
a morphism $\Phi:\Omega\to Y$ where $\Omega\subset\Tot(M)$ is open and contains the entire zero section of $M$,
\end{enumerate}
such that \(\Phi(0_{c'})=b(\rho(c'))\) for every \(c'\in C'\) and, for general $c'\in C'$, the induced rational map \(M_{c'}\dashrightarrow Y\) from the fiber of $M$ over $c'$ is nonconstant. The above ingredients are summarized in the following diagram:
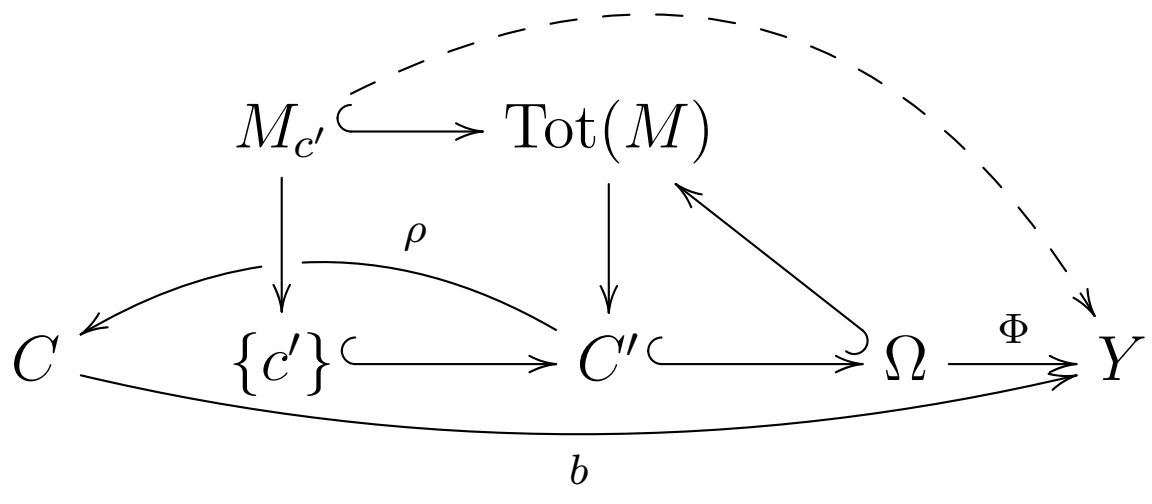
\begin{equation*}
\xymatrix{
 &M_{c'}\ar[d]\ar@{^{(}->}[r]\ar@{-->}@/^3.55pc/[rrrd]&{\Tot}(M)\ar[d] & & \\
C\ar@/_1.10pc/[rrrr]_{b} &\{c'\}\ar@{^{(}->}[r]&C'\ar@{^{(}->}[r]\ar@/_1.60pc/[ll]_(0.35){\rho}|(0.57){\hole} & \Omega\ar@{^{(}->}[ul]\ar[r]^{\Phi} & Y
}
\end{equation*}
\end{externaltheorem}

We briefly explain this result. For a fixed projective embedding of \(Y\), put 
\(a=\deg b^*\OO_Y(1)\) and let \(\delta>0\) bound the degrees of its defining equations. The jet order \(k\) in \cite[Theorem~4.2]{GX} must satisfy \(d\sum_{\ell=1}^k1/\ell>2(\dim Y+1)\delta a\). Such a \(k\) exists since \(d>0\). \cite[Proposition 3.2]{GX} implies \(\deg M>0\). We use that positivity to keep a whole affine fiber inside the prescribed open target.

We need the following consequence of the above result.
\begin{lemma}\label{lem:open-line}
Let \(U\) be a smooth quasi-projective variety, \(C\) a smooth connected projective curve, and 
\(b:C\to U\) a morphism with \(\deg b^*T_U>0\). Then there is a nonempty open \(C^0\subset C\) such that for every $\mathbb{C}$-point \(c\) of \(C^0\), there is a nonconstant morphism 
\(h:\mathbb{A}^1\to U\) with \(h(0)=b(c)\).
\end{lemma}

\begin{proof}
Choose a smooth projective compactification \(Y\) by resolving a projective closure of \(U\) away from \(U\). Since \(T_Y|_U=T_U\), the positivity assumption is precisely the degree condition for
\(C\to U\hookrightarrow Y\) in the above result. Apply the above result to obtain $\rho:C'\to C$, $M\to C'$, and $\text{Tot}(M)\supset \Omega\xrightarrow{\ \Phi \ } Y$. Let \(\Omega'=\Phi^{-1}(U)\subset\Tot(M)\). This open subset contains the zero section because 
\(\Phi(0_{c'})=b(\rho(c'))\).

To determine what has been removed from each fiber, we compactify \(\Tot(M)\) to  $$S=\PP_{\mathrm{lines}}(\OO_{C'}\oplus M)\xrightarrow{\ \pi \ } C'.$$
The two summands define the zero and infinity sections \(S_0\) and \(S_\infty\), respectively. Projection of the tautological line onto \(M\) vanishes on \(S_0\), and projection onto 
\(\OO_{C'}\) vanishes on \(S_\infty\). Taking the ratio of their divisor line bundles gives
\[
\OO_S(S_0-S_\infty)\simeq\pi^*M .
\]
If a curve in \(\Tot(M)\setminus\Omega'\) dominated \(C'\), its closure \(\Gamma\) in \(S\) would be disjoint from \(S_0\), since all points added in taking the closure lie on \(S_\infty\). But 
\(\Gamma\) is distinct from \(S_\infty\), so \(\Gamma\cdot S_\infty\ge0\). The above equation then gives
\[
  0=\Gamma\cdot S_0=\Gamma\cdot S_\infty+\deg(\Gamma\to C')\deg M>0.
\]
This is a contradiction. Thus the complement \(\Tot(M)\setminus\Omega'\) does not contain curves dominanting $C'$; it lies over finitely many points of \(C'\). Outside these points and the finitely many points where fiberwise nonconstancy of $\Phi$ fails, \(\Phi\) maps the entire affine fiber \(M_{c'}\) nonconstantly into \(U\). Choose a linear coordinate on this fiber with origin 
\(0_{c'}\). Removing the images of the exceptional points under the finite map \(\rho:C'\to C \) gives the required \(C^0\).
\end{proof}

\section{Preparation on root stacks}
We will use a root stack in two ways: its canonical divisor includes the coefficients of the orbifold boundary divisor, and a morphism to it from an ordinary curve imposes divisibility on the contact orders. 

\subsection{Morphisms from curves to root stacks}\label{sec:roots}
Let \(Z\) be a smooth projective variety, \(B=\sum_{j=1}^h B_j\) a reduced simple normal crossings divisor, and \(n_j\ge2\) integers. Let
\[
 p:\cZ=\sqrt[n_1]{(Z,B_1)}\mathbin{\times_Z}\cdots
        \mathbin{\times_Z}\sqrt[n_h]{(Z,B_h)}\longrightarrow Z .
\]
For a morphism \(v:T\to Z\), a lift to \(\cZ\) consists of a line bundle with section for each 
\(j\), together with an identification of its \(n_j\)-th tensor power with
\((v^*\OO_Z(B_j),v^*s_{B_j})\). Here \(s_{B_j}\) denotes the canonical section with divisor 
\(B_j\). For each $j$, we write \(\mathcal R_j\) for the reduced root divisor on \(\cZ\) and 
\(R_j=\OO_{\cZ}(\mathcal R_j)\) for the corresponding $n_j$-th root line bundle.

We refer to \cite[Section~2]{Cadman} and \cite[Appendix B]{AGV} for detailed treatments of the root construction.

At a point of \(Z\), take \'etale coordinates in which the components $B_i$ through the point are given by \(x_j=0\). We then have a chart for \(\cZ\) obtained by adjoining coordinates \(u_j\) with \(x_j=u_j^{n_j}\) and taking the quotient by the product of the cyclic groups \(\mu_{n_j}\), each acting on its own \(u_j\). In particular, $\cZ$ is a smooth proper Deligne--Mumford stack, and $p$ is an isomorphism over \(Z\setminus B\). By the equations \(dx_j=n_j u_j^{n_j-1}du_j\) and
\(dx_j/x_j=n_j\,du_j/u_j\), we obtain
\begin{equation}\label{eq:rootcanonical}
 p^*B_j=n_j\mathcal R_j,\quad  K_{\cZ}\sim_{\mathbb{Q}}p^*\!\left(K_Z+\sum_j(1-1/n_j)B_j\right),
\end{equation}
\begin{equation}\label{eq:kummerlog}
 p^*T_Z(-\log B)\simeq T_{\cZ}\bigl(-\log\textstyle\sum_j\mathcal R_j\bigr).
\end{equation}       
It follows that 
$$K_{\cZ}+\sum_j\mathcal R_j\sim_{\mathbb{Q}}p^*(K_Z+B).$$ 

Suppose now that \(T\) is a smooth scheme curve and, in the map $v:T\to Z$ above, the generic point of \(T\) maps outside \(B\). A lift of $v$ to $T\to \cZ$ exists if and only if every coefficient of \(v^*B_j\) is divisible by \(n_j\), for each \(j\). Indeed, taking the divisor of a root section gives the divisibility. Conversely, if \(v^*B_j=n_j A_j\) for an effective integral divisor 
\(A_j\), the pair \((\OO_T(A_j),s_{A_j})\) gives the required root. We do this for each component $B_j$, even when a point of \(T\) maps to a intersection of two divisor components. In particular, any positive contact of the coarse map of an scheme curve to \(\cZ\) is at least the corresponding 
\(n_j\).

\subsection{Incidence on a root stack}\label{sec:incidence}

We pass from Lemma~\ref{lem:open-line} to a root stack by choosing frames for its root line bundles. A root line bundle on the given projective curve can have nonzero degree. We therefore take frames relative to line bundles pulled back from an abelian variety \(\mathsf A\). The first step is to obtain these bundles after a finite cover of the curve.

\begin{lemma}\label{lem:elliptic-pencils}
Let \(M_1,\ldots,M_r\) be line bundles on a smooth connected projective curve \(C\). Then there exist a finite surjective morphism \(p:C'\to C\) from a smooth connected projective curve $C'$, a morphism \(a:C'\to \mathsf A\) to an
abelian variety,
$$C\xleftarrow{\ p\ }C'\xrightarrow{\ a\ }\mathsf A.$$
and line bundles \(L_1,...,L_r\) on \(\mathsf A\) such that \(p^*M_i\simeq a^*L_i\) for all \(i\).
\end{lemma}

\begin{proof}
Choose a line bundle \(J\) on $C$ of sufficiently large degree, so that \(N_0=J\) and \(N_i=J\otimes M_i\), \(1\le i\le r\), are very ample. For each \(j\), take a nonzero section of \(N_j\) and a second section nonvanishing on the zero set of the first. The latter choice is possible by global generation, since this zero set is finite. The two sections
give a morphism 
$$\psi_j:C\to\PP^1,\quad \text{with}\quad N_j\simeq\psi_j^*\OO_{\PP^1}(1).$$

Fix a double cover \(e:\mathsf E\to\PP^1\) by an elliptic curve $\mathsf E$. Let 
$$\mathsf A=\mathsf E^{r+1}\xrightarrow{\ e^{\times r+1}\ }(\mathbb{P}^1)^{r+1}.$$ 
Consider the fiber product $C\times_{(\mathbb{P}^1)^{r+1}}\mathsf{A}$ via the maps $e^{\times r+1}$ and $(\psi_0,...,\psi_r)$. Normalize the closure of a point of the generic fiber, with its reduced structure. The resulting curve $C'$ is smooth, connected and projective, and its projections 
\(p:C'\to C\) and \(a:C'\to \mathsf A\) satisfy $e\mathbin{\!\scriptstyle\circ\!}\mathrm{pr}_j\mathbin{\!\scriptstyle\circ\!} a=\psi_j\mathbin{\!\scriptstyle\circ\!} p.$ The ingredients above are summarized in the following diagram:
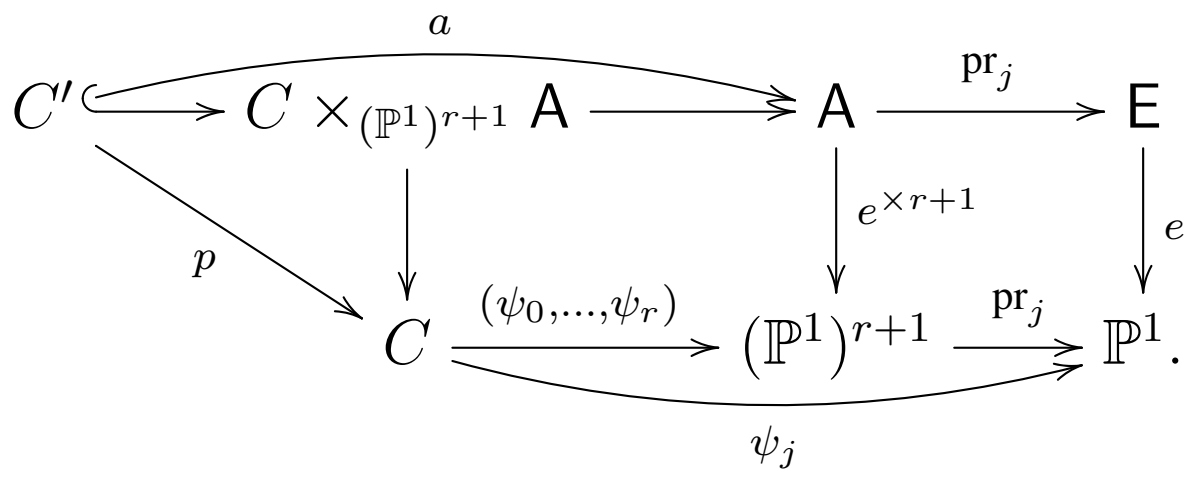
\begin{equation*}
\xymatrix{
C'\ar@{^{(}->}[r]\ar[dr]_p\ar@/^0.90pc/[rr]^a & C\times_{(\mathbb{P}^1)^{r+1}}\mathsf{A}\ar[d]\ar[r] &\mathsf{A}\ar[d]^(0.4){e^{\times r+1}}\ar[r]^{\text{pr}_j}& \mathsf{E}\ar[d]^{e} \\
 &C\ar[r]^(0.4){(\psi_0,\ldots,\psi_r)}\ar@/_0.90pc/[rr]_{\psi_j} & (\mathbb{P}^1)^{r+1}\ar[r]^(0.6){\text{pr}_j}&\PP^1. 
}    
\end{equation*}
In particular, \(p\) is finite and surjective. If \(P_j=\mathrm{pr}_j^*e^*\OO_{\PP^1}(1)\), then 
\(a^*P_j\simeq p^*N_j\). Taking \(L_i=P_i\otimes P_0^{-1}\) gives the assertion.
\end{proof}

\begin{proposition}\label{prop:incidence}
Let $$\cZ\to Z$$ 
be the root stack of a smooth projective variety $Z$ along a reduced SNC divisor \(B=\sum_{i=1}^r B_i\) with orders \(\{n_i\}\).
Let \(f:C\to\cZ\) be a morphism from a smooth connected projective curve whose coarse image meets \(Z\setminus B\). Suppose that \(-\deg f^*K_{\cZ}>0\). Then there is a nonempty open subset $$C^0\subset C,$$
such that, for every $\mathbb{C}$-point \(c\) of \(C^0\), there is a morphism \(g:\PP^1\to\cZ\) with nonconstant
coarse map and \(g(0)\simeq f(c)\).
\end{proposition}

\begin{proof}
We divide the proof into several parts.

\noindent\emph{Construction of the scheme of frames.}
Let \(R_i\) be the \(i\)-th root line bundle on \(\cZ\). So $R_i^{\otimes n_i}$ is isomorphic to the pull-back of $\mathcal{O}_Z(B_i)$. We apply Lemma~\ref{lem:elliptic-pencils} to \(M_i=f^*R_i, i=1,...,r\), and choose
isomorphisms 
\begin{equation}\label{eqn:isom_chosen}
a^*L_i\simeq p^*f^*R_i, \quad i=1,...,r.
\end{equation}
 Over \(\mathsf A\times\cZ\), form the
fiber product of the Isom spaces:
\begin{equation}\label{eqn:torsor_frames}
 W=\prod_{i=1}^r \mathrm{Isom}_{\mathsf A\times \cZ}\bigl(\mathrm{pr}_{\mathsf A}^*L_i,\mathrm{pr}_{\cZ}^*R_i\bigr)\to \mathsf A\times\cZ.
\end{equation}
The structure map of (\ref{eqn:torsor_frames}) is a \((\mathbb{G}_m)^r\)-torsor, with the factors of the torus scaling the respective frames. The chosen isomorphisms (\ref{eqn:isom_chosen}) give a lift \(\widetilde f:C'\to W\) of 
\((a,f\mathbin{\!\scriptstyle\circ\!} p): C'\to \mathsf{A}\times \cZ\):
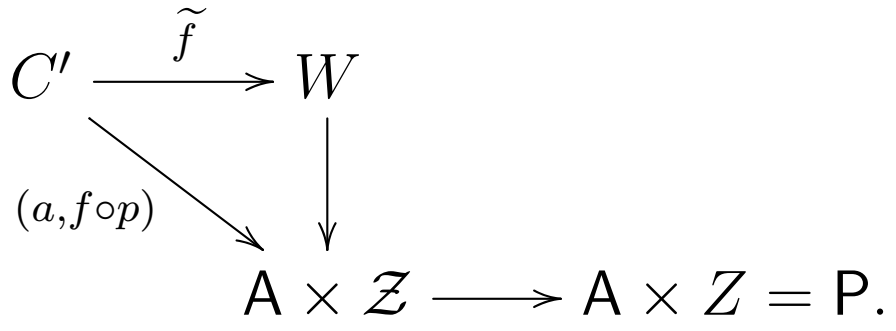
\begin{equation*}
\xymatrix{
C'\ar[r]^{\widetilde{f}}\ar[dr]_(0.4){(a,f\mathbin{\scriptstyle\circ} p)} &W\ar[d] & \\
 & \mathsf{A}\times \cZ\ar[r] & \mathsf{A}\times Z=\mathsf P. 
}    
\end{equation*}

We verify that \(W\) is a quasi-projective scheme by writing it over \(\mathsf P=\mathsf{A}\times Z\). On \(\mathsf P\), put
\[
 \Lambda_i=\mathrm{pr}_A^*L_i,\qquad H_i=\mathrm{pr}_Z^*\OO_Z(B_i),
 \qquad \mathcal H_i=\mathcal{H}\!om_{\mathsf P}(\Lambda_i^{\otimes n_i},H_i),
\]
and let \(\sigma_i\) be the pullback of the canonical section of \(\OO_Z(B_i)\). Set 
$$E_{\mathsf P}=\bigoplus_{i=1}^r(\mathcal H_i\oplus\Lambda_i).$$
A point of \(\Tot(E_{\mathsf P})\) over a scheme \(T\to \mathsf P\) consists of homomorphisms 
\(\eta_i:\Lambda_{i,T}^{\otimes n_i}\to H_{i,T}\) and sections \(v_i\) of \(\Lambda_{i,T}\), where a subscript \(T\) denotes pullback. Impose the following conditions:
\[
\eta_i\text{ is invertible},
\qquad
\eta_i(v_i^{\otimes n_i})=\sigma_{i,T}, \quad 1\le i\le r.
\]
The first conditions are open. On the open set defined by these conditions, the second conditions cut out a closed subscheme \(W'\), since each is the vanishing of a section of a line bundle. Thus \(W'\) is locally closed in \(\Tot(E_{\mathsf P})\). Since $\Tot(E_{\mathsf P})$ can be compactified to the projective bundle of lines \(\PP_{\mathsf{P}}(\OO_{\mathsf{P}}\oplus E_{\mathsf{P}})\), 
\(W'\) is quasi-projective.

To identify \(W'\) with \(W\), write a framed root over \(T\) as 
$$(M_i,\tau_i,\phi_i,\theta_i),$$ 
where \(\phi_i:M_i^{\otimes n_i}\xrightarrow{\ \sim\ }H_{i,T}\) mpas \(\tau_i^{\otimes n_i}\) to \(\sigma_{i,T}\), and \(\theta_i:\Lambda_{i,T}\xrightarrow{\ \sim\ }M_i\) is the frame. Its coordinates in \(W'\) are
\[
 v_i=\theta_i^{-1}(\tau_i),\quad \eta_i=\phi_i\mathbin{\!\scriptstyle\circ\!}\theta_i^{\otimes n_i}.
\]
In the other direction, we take the root \((\Lambda_{i,T},v_i,\eta_i)\) with the identity frame. The isomorphism from this root to the original framed root must be \(\theta_i\), including at a zero of \(\tau_i\). These assignments
commute with pullback and are inverse equivalences on every test scheme. This proves \(W\simeq W'\).

The equations above also describe the map which forgets the \(v_i\). Let \(V\subset\Tot(\bigoplus_i\mathcal H_i)\) be the locus where all \(\eta_i\) are invertible. In local trivializations on \(\operatorname{Spec}R\subset V\), write 
\(\eta_i=\lambda_i\in R^\times\) and \(\sigma_i=b_i\in R\). The algebra of \(W\) is then
\[
 R[t_1,\ldots,t_r]\big/\bigl(t_i^{n_i}-\lambda_i^{-1}b_i \mid 1\le i\le r\bigr).
\]
The monomials \(\prod_i t_i^{e_i}\), \(0\le e_i<n_i\), form a basis over \(R\), so \(W\to V\) is  locally free of rank \(\prod_i n_i\). Smoothness of \(W\) follows from its torsor map to the smooth stack \(\mathsf A\times\cZ\).

\noindent\emph{An affine line with nonconstant coarse projection.}
The tangent bundle along the fibers of the torsor (\ref{eqn:torsor_frames}) is trivial, with basis given by the infinitesimal actions of the torus factors. Pulling back its tangent sequence to $C'$ gives
\[
 0\longrightarrow\OO_{C'}^{\oplus r}
   \longrightarrow\widetilde f^*T_W
   \longrightarrow a^*T_{\mathsf A}\oplus(f\mathbin{\!\scriptstyle\circ\!} p)^*T_{\cZ}
   \longrightarrow0 .
\]
\(T_{\mathsf A}\) can be trivialized, hence \(\deg\widetilde f^*T_W=(\deg p)\cdot(-\deg f^*K_{\cZ})>0\). Applying 
Lemma~\ref{lem:open-line} to the connected component of \(W\) containing \(\widetilde f(C')\), we obtain a nonempty open \(U'\subset C'\) such that, for each \(c'\in U'(\mathbb{C})\), there is a nonconstant map 
$$h:\mathbb A^1\to W$$ 
with \(h(0)=\widetilde f(c')\). Shrink \(U'\) so that \(f(p(c'))\in \cZ\) lies over \(Z\setminus B\).

The composition of \(h\) and $\text{pr}_{\mathsf{A}}: W\to\mathsf A$ extends to a map $\PP^1\to \mathsf{A}$ since $\mathsf{A}$ is proper. This map is constant because an abelian variety contains no rational curves. If the composition of \(h\) and $W\to Z$ were also constant, \(h\) would factor through a fiber of \(W\to P=\mathsf{A}\times Z\) over 
\(\mathsf A\times(Z\setminus B)\). The root stack is isomorphic to the coarse space there, so this fiber, with base point \(h(0)\), is a copy of \((\mathbb{G}_m)^r\). Thus \(h\) would be constant, a contradiction. Therefore $h$ is not constant after projecting to $Z$.

\noindent\emph{Compactification of the source.}
Because $Z$ is proper, the projection of \(h\) to $Z$ extends to a noncontant map 
$$F:\PP^1\to Z.$$
Since \(F(0)\notin B\), all \(F^*B_i\) are effective divisors, and their coefficients at finite points are divisible by \(n_i\) by the existing lift on \(\mathbb A^1\). Let \(\ell\) be a common multiple of the \(n_i\), taking \(\ell=1\) in case $B=\emptyset$. The map \(q:\PP^1\to\PP^1\), 
\(t\mapsto t^\ell\), fixes \(0\) and is totally ramified at infinity. Thus every finite multiplicity in \((F\mathbin{\!\scriptstyle\circ\!} q)^*B_i\) is still divisible by \(n_i\), and the multiplicity at infinity is \(\ell\operatorname{ord}_\infty F^*B_i\). The lifting criterion of Section~\ref{sec:roots} gives $$g:\PP^1\to\cZ$$
with coarse map \(F\mathbin{\!\scriptstyle\circ\!} q\). This map is nonconstant, and its value at \(0\) is isomorphic to 
\(f(p(c'))\), since that point lies over \(Z\setminus B\). 
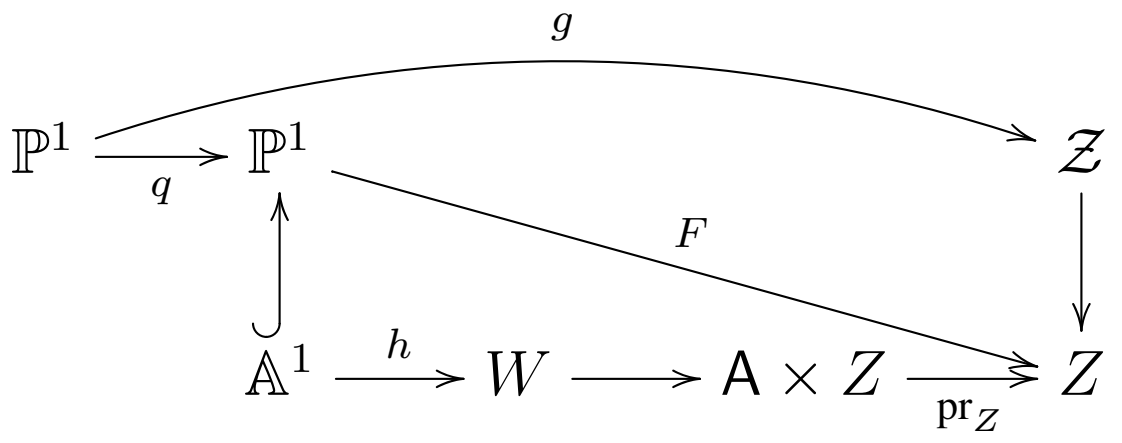
\begin{equation*}
\xymatrix{
\PP^1\ar[r]_{q}\ar@/^1.50pc/[rrrr]^{g} & \PP^1\ar[drrr]^{F} &{}&{}&\cZ\ar[d]\\
{}&\mathbb{A}^1\ar@{^{(}->}[u]\ar[r]^{h} & W\ar[r] &\mathsf{A}\times Z\ar[r]_(0.6){\text{pr}_Z}& Z 
}    
\end{equation*}
Finally, take 
$$C^0=C\mathbin{\big\backslash} p(C'\setminus U').$$ 
The omitted set is finite and \(p\) is surjective, so this choice works for every 
\(c\in C^0(\mathbb C)\).
\end{proof}
\begin{remark}
By applying bend-and-break argument to the stack $\cZ$ as in e.g. \cite{CT_cone}, we can get a nonconstant map $\mathsf{f}:\mathcal{C}\to\cZ$ from a smooth genus $0$ {\em orbifold} curve $\mathcal{C}$. Taking coarse moduli spaces gives a nonconstant map $\bar{\mathsf{f}}:\PP^1\to Z$. It is possible to control the number of orbifold points on $\mathcal{C}$, but we were not able to gaurantee that there is no orbifolds on $\mathcal{C}$.      
\end{remark}

\section{Testing an effective divisor with a free curve}\label{sec:testing}

Fix a nonzero effective Cartier divisor \(E\) on \(X\). The goal here is to apply Proposition~\ref{prop:incidence} on an auxiliary space over \(X\), chosen so that the resulting rational curves have positive \(E\)-degree. 

\subsection{A covering family}
We first explain how incidence on a complete curve gives a covering family.
\begin{lemma}\label{lem:covering}
Let $$p:\mathcal S\to S$$
be the root stack of a smooth projective variety $S$ along an SNC divisor with finite root indices. Write
$$K_{\mathcal S}\sim_{\mathbb Q}p^*M.$$
Suppose that \(M\cdot\beta<0\) for a complete-intersection curve class \(\beta\) cut out by sufficiently ample divisors on \(S\). In dimension one we take \(\beta=[S]\). Then there is a dominating family of morphisms 
$$g:\PP^1\to\mathcal S$$
with nonconstant coarse maps. Its general member is \emph{free on the stack}, meaning that \(g^*T_{\mathcal S}\) is globally generated, and satisfies \(-\deg g^*K_{\mathcal S}\ge2\).
\end{lemma}

\begin{proof}
Fix a projective embedding of \(S\). The root stack $\mathcal{S}$ is tame over \(\mathbb C\) and has projective coarse space \(S\). For each positive coarse degree, the stack of unmarked genus-$0$ twisted stable maps to \(\mathcal S\)
is of finite type by \cite[Theorem~1.4.1]{AV}. On the locus parametrizing maps with smooth domains, the domains have neither nodes nor markings, so the domains are ordinary genus-$0$ curves. Choices of parametrizations of this curve form a \(\operatorname{PGL}_2\)-torsor. Taking scheme atlases and their reduced irreducible components gives
countably many families of parametrized maps, indexed over all positive degrees. They include every 
\(g:\PP^1\to\mathcal S\) with nonconstant coarse map, since such a map is stable.

Suppose that the closures \(Z_j\subset S\) of the coarse evaluation images were all proper subschemes. Choose a $\mathbb{C}$-point \(x\) outside the boundary divisor and all the \(Z_j\), using the uncountability of \(\mathbb C\). By Bertini, a general complete-intersection curve \(B\subset S\) of class \(\beta\) through \(x\) is smooth and irreducible, meets the boundary divisor transversely, and avoids its deeper strata. Here the defining linear systems through \(x\) have no other base point and separate tangent directions at \(x\). In dimension $1$, we take \(B=S\). Since \(x\notin Z_j\), each intersection \(B\cap Z_j\) is finite.

To apply Proposition~\ref{prop:incidence}, we lift a finite cover of \(B\) to \(\mathcal S\). Let 
\(P\) be the finite set of intersections with boundary divisor. If \(P\ne\emptyset\), let \(N\) be a common multiple of the root indices. Choose a rational function \(u\) on \(B\) with a simple zero at every point of \(P\). For example, for \(z\in B\setminus P\) and \(k\) sufficiently large, Riemann--Roch implies that we can prescribe the first jets of sections of \(\OO_B(kz)\) at all points of \(P\). Normalize \(B\) in the extension defined by \(t^N=u\). At each point of \(P\) this equation is Eisenstein, so the resulting smooth connected curve \(B'\) has degree \(N\) over \(B\) and is totally ramified there. If \(P=\emptyset\), we take \(B'=B\). In either case, the pullback of each divisor component is divisible by its root index. The lift criterion in Section \ref{sec:roots} gives \(b:B'\to\mathcal S\) with
\[
\deg b^*K_{\mathcal S}=\deg(B'/B)\,(M\cdot\beta)<0.
\]
We then apply Proposition \ref{prop:incidence} to $b: B'\to \mathcal{S}$. The coarse images of the resulting rational curves on \(\mathcal S\) contain all but finitely many points of \(B\). On the other hand, \(\bigcup_j(B\cap Z_j)\) is countable. This proves that some family dominates.

We now argue that any such family has free general member. Adjoining reparametrizations makes evaluation at \(0\in\PP^1\) dominant. Restrict to the smooth locus of the reduced parameter variety and to maps taking \(0\) outside the boundary. At a general parameter, the differential of this evaluation is surjective by generic smoothness. A tangent vector to the parameter variety induces a section of \(g^*T_{\mathcal S}\), and its image under the evaluation differential is the value of that section at \(0\). Thus global sections span the fiber there. Writing
\(g^*T_{\mathcal S}=\bigoplus_j\OO_{\PP^1}(a_j)\), we obtain \(a_j\ge0\) for every \(j\), since a negative summand has no sections. The differential \(dg:T_{\PP^1}=\OO_{\PP^1}(2)\to g^*T_{\mathcal S}\) is generically nonzero because the coarse map is nonconstant in characteristic zero. Some 
\(a_j\) is therefore at least $2$, and
\begin{equation}\label{eq:testing-degree}
-\deg g^*K_{\mathcal S}=\sum_j a_j\ge2.
\end{equation}
Both conclusions are invariant under reparametrization.
\end{proof}

\subsection{A klt perturbation on the base}
Put \(A=-(K_X+\Delta)\). We need an effective representative of a multiple of \(A-\varepsilon E\) with controlled singularities.

\begin{lemma}\label{lem:perturbation}
Suppose that \(A\) is nef and big. Then there is an effective $\mathbb{Q}$-divisor \(T\) such that 
\((X,\Delta+T)\) is klt log Fano. For every nonzero effective Cartier divisor \(E\), there also exist a rational number \(\varepsilon>0\), an integer \(m\ge3\), and an effective integral divisor 
\(G\sim_{\mathbb Q}m(A-\varepsilon E)\) such that \((X,\Delta+G/m)\) is klt.
\end{lemma}

\begin{proof}
By Kodaira's lemma (see e.g. \cite[Corollary~2.2.7]{Laz}), write 
$$A\sim_{\mathbb Q}H+G_0$$
with \(H\) ample and \(G_0\ge0\) a $\mathbb Q$-divisor. Fix a log resolution of 
\(D+\operatorname{Supp}G_0\). The finitely many coefficients in the log pullback of \(\Delta+tG_0\) depend linearly on \(t\). They are all less than $1$ at \(t=0\), since \((X,\Delta)\) is klt. Choose a rational number \(t\in (0,1)\) so that these inequalities hold, and set \(T=tG_0\). Then
\[
 -(K_X+\Delta+T)=A-T\sim_{\mathbb Q}(1-t)A+tH
\]
is ample, because it is the sum of a nef class and an ample class.

Choose a sufficiently small positive rational number \(\varepsilon\) so that 
\(P=A-T-\varepsilon E\) is ample. Choose \(m\ge3\) sufficiently divisible, so that \(mT\) is integral and \(mP\) is very ample and Cartier. For a general \(G_1\in |mP|\), put $$G=mT+G_1.$$
On the chosen resolution, the pulled-back linear system is base point free. By Bertini, its general member is smooth, contains no exceptional component, and meets the strata of the resolved boundary divisor transversely. Adding it with coefficient \(1/m<1\) therefore gives a log pullback with SNC support and coefficients less than $1$. Since \(\Delta+G/m=\Delta+T+G_1/m\), this proves the klt assertion, and \(G\sim_{\mathbb Q}m(A-\varepsilon E)\) by construction.
\end{proof}

\subsection{The hyperplane arrangement}
We use the canonical section of $G$ to vary one hyperplane in a projective bundle. Along $G$, some of the hyperplanes meet in excess dimension. The following lemma computes the effect on discrepancies.

\begin{lemma}\label{lem:arrangement}
Let \(X\) be smooth over \(\mathbb{C}\), \(\Delta\) an effective $\mathbb{Q}$-divisor, \(G\) an effective integral divisor, and \(m\ge3\) an integer. Set \(r=m-2\), \(w=1-1/m\), and 
\(L=\OO_X(G)\), with canonical section \(s_G\). In the convention of $1$-dimensional quotients, put
\[
 q:Y=\PP_X(\OO_X^{\oplus r}\oplus L)\longrightarrow X,\qquad
 H_\infty=\PP_X(\OO_X^{\oplus r}).
\]
Let $$H_1,\ldots,H_{r+1}$$ 
be the zero divisors of the images in the tautological quotient of 
$$e_1,\ldots,e_r, \text{ and } e_1+\cdots+e_r+s_G,$$ where \(e_j\) are the standard basis sections. If \((X,\Delta+G/m)\) is klt, then so is \((Y,B_Y)\), where
\[
B_Y=q^*\Delta+w\left(H_\infty+\sum_{j=1}^{r+1}H_j\right).
\]
\end{lemma}

\begin{proof}
We compare \((Y,B_Y)\) with a constant arrangement by a birational change of the fiber coordinates. We work over an open set \(U\subset X\) on which \(L\) is trivial, \(s_G\) is given by a regular function \(s\), and \(K_X\) has a nowhere vanishing generator \(\eta\). Write 
$$Y_1=Y|_U, B_1=B_Y|_{Y_1}, \text{ and } q_1:Y_1\to U.$$
On the complement of \(H_\infty\), use affine coordinates \(y_1,\ldots,y_r\). Consider the rational function $P_y$ and top form $\alpha_y$ where
\[
P_y=y_1\cdots y_r\bigl(\textstyle\sum_j y_j+s\bigr),
\quad
\alpha_y=\eta\wedge dy_1\wedge\cdots\wedge dy_r.
\]
The factors of \(P_y\) vanish on the \(H_j\) and have simple poles along \(H_\infty\). Each stays a nonzero polynomial on every fiber, so there are no vertical components in its divisor. The affine volume form has pole order \(r+1\) at infinity, as is seen by putting \(y_1=1/z\), \(y_j=u_j/z\) for \(j>1\). Thus
\[
\operatorname{div}(P_y)=\sum_{j=1}^{r+1}H_j-(r+1)H_\infty,
\qquad \operatorname{div}(\alpha_y)=-(r+1)H_\infty .
\]
Taking \(K_{Y_1}=\operatorname{div}(\alpha_y)\) and using \((r+2)w=r+1\), we obtain the equality of rational divisors
\[
K_{Y_1}+B_1=q_1^*\Delta+w\operatorname{div}_{Y_1}(P_y).
\]

On \(Y_2=U\times\PP^r\), let \(\mathcal H\) be the constant arrangement defined in homogeneous coordinates \([t:v_1:\cdots:v_r]\) by \(t=0\), \(v_j=0\), and \(t+\sum_jv_j=0\). Write \(q_2:Y_2\to U\) and 
$$B_2=q_2^*(\Delta+G/m)+w\mathcal H.$$
On the chart \(t=1\), put \(P_v=v_1\cdots v_r(\sum_jv_j+1)\) and 
\(\alpha_v=\eta\wedge dv_1\wedge\cdots\wedge dv_r\). The same calculation gives
\[
\operatorname{div}_{Y_2}(\alpha_v)+B_2 =q_2^*(\Delta+G/m)+w\operatorname{div}_{Y_2}(P_v).
\]
The substitution \(y_j=sv_j\) identifies the two function fields over \(U\). Under this identification, \(P_y=s^{r+1}P_v\) and \(\alpha_y=s^r\alpha_v\), since \(\eta\wedge ds=0\).
Using the same top form \(\alpha_y\) for \(K_{Y_2}\), we find
\[
 \begin{aligned}
 K_{Y_2}+B_2
 &=q_2^*\Delta+(r+1/m)q_2^*G+w\operatorname{div}_{Y_2}(P_v)\\
 &=q_2^*\Delta+w\operatorname{div}_{Y_2}(P_y),
 \end{aligned}
\]
where the second equality uses \((r+1)w=r+1/m\).

Take a common log resolution of the two pairs $(Y_1,B_1)$ and $(Y_2, B_2)$ by resolving the closure of the graph $\{y_j=sv_j, j=1,...,y\}$ over \(U\). The projections from the graph are proper because both models are projective over \(U\). The expressions for \(K_{Y_i}+B_i\) have equal pullbacks to the resolution because the two maps to \(U\) agree there and \(P_y\) is the same rational function. Thus all discrepancies agree. On \(Y_2\), any \(r+1\) of the \(r+2\) homogeneous defining forms are independent, so \(\mathcal H\) has SNC support. The product of a log resolution of \((U,(\Delta+G/m)|_U)\) with \(\PP^r\) resolves \((Y_2,B_2)\). The log pullback has the coefficients of the resolved base pair and the coefficient \(w<1\) on the fiber hyperplanes. All are less than $1$. This proves that \((Y_2,B_2)\), and hence \((Y_1,B_1)\), is klt. Repeating the argument to this kind of open sets that cover \(X\) proves the lemma.
\end{proof}

\subsection{Positive degree against a prescribed divisor}
We can now proceed to establish the main result of this section. 

Let $X, D, \Delta$ be as in Section \ref{sec:result}. Let $\pi:\cX\to X$ be the stack obtained by taking the \(m_i\)-th root along each \(D_i\).

\begin{proposition}\label{prop:testing}
Suppose that \(A=-(K_X+\Delta)\) is nef and big. For a nonzero effective Cartier divisor\footnote{\(E\) can be one of the \(D_i\).} \(E\) on \(X\), there is a dominating family of morphisms 
\(f:\PP^1\to\cX\) whose general member is free on the stack and satisfies \(E\cdot(\pi\mathbin{\!\scriptstyle\circ\!} f)_*[\PP^1]>0\).
\end{proposition}

\begin{proof}
Let \(\varepsilon,m,G\) be as in Lemma~\ref{lem:perturbation} and form the klt pair \((Y,B_Y)\) as in Lemma~\ref{lem:arrangement}. We first compute its log canonical class. Put 
\(\xi=c_1(\OO_Y(1))\). The map from \(q^*L\) to the tautological quotient cuts out \(H_\infty\), of class \(\xi-q^*G\). Each of the other \(r+1\) defining sections is nowhere zero in 
\(\mathcal{O}_X^{\oplus r}\oplus L\), so its zero divisor in the tautological quotient is a projective subbundle of class \(\xi\). The \(m=r+2\) hyperplanes are distinct over 
\(X\setminus\operatorname{Supp}G\). Their total class is \(m\xi-q^*G\), while 
$$K_Y\sim q^*(K_X+G)-(m-1)\xi.$$
Since \(mw=m-1\), adding \(B_Y\) gives
\begin{equation}\label{eq:auxcan}
K_Y+B_Y\sim_{\mathbb Q}q^*(K_X+\Delta+G/m)\sim_{\mathbb Q}-\varepsilon q^*E .
\end{equation}

Consider a projective log resolution of \((Y,B_Y)\): 
$$\mu:\widetilde Y\to Y.$$
For each original prime component \(D_i\) of \(D\subset X\), the inverse image \(q^*D_i=Y|_{D_i}\) is a reduced prime divisor, because it is a projective bundle over \(D_i\). Denote its strict transform to $\widetilde{Y}$ by \(D_i'\). Let \(B'\) be the strict transform of \(B_Y\) to $\widetilde{Y}$ with its coefficients and let \(F_j\subset \widetilde{Y}\) be the exceptional prime divisors of $\mu$. With compatible canonical divisors, write
\[
 K_{\widetilde Y}+B'=\mu^*(K_Y+B_Y)+\sum_j a_jF_j.
\]
We have \(a_j>-1\) for every $j$ because $(Y, B_Y)$ is klt. Choose \(N\ge2\) divisible by all the 
\(m_i\) and sufficiently large, so  that \(a_j+1-1/N\ge0\) for each \(j\). On \(\widetilde Y\), the strict transforms and the \(F_j\)'s have SNC support. Form the root stack 
$$p:\mathcal S\to\widetilde Y$$ 
with order \(m_i\) along \(D_i'\), order \(m\) along each strict transform of a hyperplane, and order \(N\) along each \(F_j\). By (\ref{eq:rootcanonical}) and (\ref{eq:auxcan}), its canonical class is $$K_{\mathcal S}\sim_{\mathbb Q}p^*M,$$ where,
\begin{equation}\label{eq:exceptional}
\begin{aligned}
M=K_{\widetilde Y}+B'+(1-1/N)\sum_j F_j
      \sim_{\mathbb Q}-\varepsilon\mu^*q^*E+R, \quad 
R=\sum_j(\underbrace{a_j+1-1/N}_{\geq 0 \text{ by choice of }N})F_j .
 \end{aligned}
\end{equation}

The choice of \(N\) also gives a morphism $\widetilde q:\mathcal S\to\cX$ over 
\(Q=q\mathbin{\!\scriptstyle\circ\!}\mu\mathbin{\!\scriptstyle\circ\!}  p:\mathcal S\to X\),
\begin{equation*}
\xymatrix{
 & & & \mathcal{X}\ar[d]^{\pi}\\
\mathcal{S}\ar@/^0.90pc/[urrr]^{\widetilde{q}}\ar[r]^{p}\ar@/_0.90pc/[rrr]_{Q} &\widetilde{Y}\ar[r]^{\mu}&Y\ar[r]^{q} & X.
}    
\end{equation*}
Indeed, \(\mu^*q^*D_i=D_i'+\sum_j b_{ij}F_j\) for integers \(b_{ij}\ge0\). Write \(\mathcal D_i,\mathcal F_j\) for the reduced root divisors on \(\mathcal S\) above \(D_i',F_j\), respectively. Since \(p^*D_i'=m_i\mathcal D_i\) and \(p^*F_j=N\mathcal F_j\), we have
\[
\frac{1}{m_i}Q^*D_i=\mathcal D_i+\sum_j b_{ij}\frac{N}{m_i}\mathcal F_j.
\]
The right-hand side is an effective integral Cartier divisor. Its line bundle with canonical section gives the required \(m_i\)-th root of \((Q^*\OO_X(D_i),Q^*s_{D_i})\). Applying the universal property for all \(i\) gives \(\widetilde q\).

It remains to use the two properties of \(R\) in \eqref{eq:exceptional}. First, \(R\) is \(\mu\)-exceptional. For \(n=\dim Y\) and a very ample divisor \(J\) on \(Y\), since \(q^*E\) is nonzero and effective, the projection formula then gives
\[
M\cdot(\mu^*J)^{n-1} =(-\varepsilon\,q^*E)\cdot J^{n-1}<0.
\]
For an ample divisor \(J'\) on \(\widetilde Y\), the class \(\mu^*J+\delta J'\) is ample for every 
\(\delta>0\). By continuity, its \((n-1)\)-st power still has negative \(M\)-degree for sufficiently small rational \(\delta>0\). A sufficiently divisible multiple of it is then very ample. We then apply Lemma~\ref{lem:covering} to obtain a dominating family $$g:\PP^1\to\mathcal S$$
with \(-\deg g^*K_{\mathcal S}\ge2\) for a general member.

Second, \(R\) is effective by our choice of \(N\). For a general member \(g\), the curve \(p\mathbin{\scriptstyle\circ\!} g\) meets the complement of \(\operatorname{Supp}R\), so \(\deg(p\mathbin{\!\scriptstyle\circ\!} g)^*R\ge0\). Taking degrees after pullback of \eqref{eq:exceptional} by \(p\mathbin{\!\scriptstyle\circ\!} g\), we obtain
\[
\varepsilon\,E\cdot(Q\mathbin{\!\scriptstyle\circ\!} g)_*[\PP^1] =-\deg g^*K_{\mathcal S}+\deg(p\mathbin{\!\scriptstyle\circ\!} g)^*R\ge2.
\]
In particular, \(f=\widetilde q\mathbin{\!\scriptstyle\circ\!} g\) has nonconstant coarse map and the required positive degree. The family obtained from this by projecting to $X$ is dominant, because \(Q\) is surjective. We apply the argument in the last paragraph of the proof of Lemma~\ref{lem:covering} to this projected family to conclude that its general member is free on \(\cX\).
\end{proof}

\section{From covering families to Campana rational connectedness}\label{sec:completion}
We complete the proof of Theorem \ref{thm:main} in this section.

\subsection{Choosing the curve class}
Lemma~\ref{lem:perturbation} gives a klt log Fano divisor on the smooth projective variety \(X\). Thus \(X\) is a Mori dream space by \cite[Corollary~1.3.2]{BCHM}. In particular, there are nonzero effective Cartier divisors 
$$E_1,\ldots,E_a$$ 
whose numerical classes generate \(\overline{\operatorname{Eff}}^1(X)\). Indeed, take the zero divisors of nonconstant homogeneous Cox ring generators. Every effective class has a section in the Cox ring, after taking a positive multiple if one uses a finite-index grading lattice. A nonzero homogeneous section has a monomial of its degree in the generators. Their degrees therefore span the effective cone, which is closed by finite generation. Each nonzero effective divisor has positive degree against an ample complete-intersection class.

For each \(j\), take an irreducible covering family as in Proposition~\ref{prop:testing} for \(E_j\), and denote its coarse curve class by \(\beta_j\). The class $\beta_j$ is constant on the family, since the degree of a line bundle is constant in a flat family of proper curves. For any effective divisor \(E\), we can choose a member of the family whose image meets \(X\setminus E\). Its pullback of \(E\) is effective, so \(E\cdot\beta_j\ge0\). For \(E=E_j\), Proposition \ref{prop:testing} gives the strict inequality. Consequently,
\[
\beta:=\sum_{j=1}^a\beta_j,
\qquad
E_k\cdot\beta\ge E_k\cdot\beta_k>0, \quad 1\le k\le a.
\]
The finitely many strict inequalities hold in a neighborhood of \(\beta\) in \(N_1(X)_{\mathbb R}\). This neighborhood lies in the dual of \(\overline{\operatorname{Eff}}^1(X)\), as required.

\subsection{Realizing the sum by a free curve}
\begin{lemma}\label{lem:smooth-sum}
Let \(a\ge1\). Suppose that 
$$h_j:\PP^1\simeq C_j\to\cX,\quad 1\le j\le a,$$
are free on the stack $\cX$ and have nonconstant coarse maps. Suppose also that there are points 
\(p_j\in C_j, j=1,...,a\) mapping to a common point \(x\in X\setminus D\). Then, there is a map $$h':\PP^1\to\cX,$$
free on the stack and passing through \(x\), whose coarse curve class is the sum of the coarse curve classes of the 
\(h_j\).
\end{lemma}

\begin{proof}
Let $C_0$ be a copy of $\PP^1$. Take distinct points \(z,q_1,\ldots,q_a\) on \(C_0\). Attach \(C_j\) to \(C_0\) by identifying \(p_j\) with \(q_j\). Define a map \(h:C\to\cX\) from the resulting nodal curve $C$ to $\cX$ by using the constant map \(C_0\to \{x\}\) and \(h_j\) on \(C_j\); since \(\pi\) is an isomorphism over \(x\), these restrictions agree at each node. Besides \(z\), mark one further point $w\in C_0$ away from the nodes. Together with a node, these markings stabilize the contracted component. See Figure \ref{fig:comb-smoothing} for a depiction.

\begin{figure}[htbp]
  \centering
  \resizebox{0.50\textwidth}{!}{%
  \begin{tikzpicture}[
    x=1cm, y=1cm,
    line cap=round,
    every node/.style={font=\small}
  ]
    \draw[thick] (-5,0) -- (5,0);
    \node[right] at (5,0) {$C_0$};

    \foreach \x/\j in {-2.8/1,-0.5/2,3.0/a} {
      \draw[thick]
        (\x-0.10,-0.35)
        .. controls (\x,-0.15) and (\x,-0.10) .. (\x,0)
        .. controls (\x,0.70) and (\x+0.15,1.40)
        .. (\x+0.55,1.80);

      \fill (\x,0) circle[radius=1.7pt];

      \node[above] at (\x+0.55,1.80) {$C_{\j}$};
      \node[below] at (\x,-0.48) {$q_{\j}\sim p_{\j}$};
    }

    \node at (1.55,1.25) {$\cdots$};

    \foreach \x/\name in {-4.3/z,4.3/w} {
      \draw[thick,fill=white]
        (\x,0) circle[radius=2.2pt];
      \node[above] at (\x,0.13) {$\name$};
    }

    \node at (0,-1.25) {$h|_{C_0}\equiv x$};
  \end{tikzpicture}%
}
\caption{\small The nodal curve $C=C_0\cup C_1\cup\cdots\cup C_a$. Each component is a copy of $\PP^1$, and $q_j\in C_0$ is identified with $p_j\in C_j$. The point $w$ is the auxiliary marking on $C_0$. The map $h$ contracts $C_0$ to $x$ and restricts to $h_j$ on $C_j$.}
  \label{fig:comb-smoothing}
\end{figure}
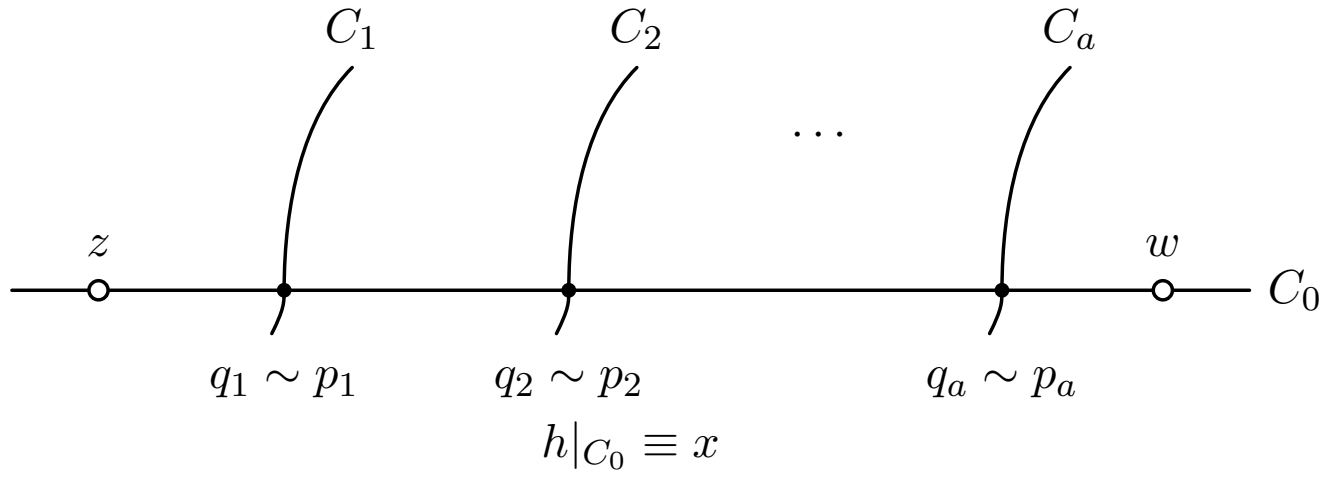

We construct a deformation of $h: C\to \cX$ that smooths \(C\) while keeping the value of $h$ at \(z\) fixed. Put 
\(V=h^*T_{\cX}\). The normalization sequence on \(C\), with the component sheaves understood to be pushed forward, reads
\[
0\longrightarrow V(-z)
 \longrightarrow
\bigl(T_xX\otimes\OO_{C_0}(-z)\bigr)\oplus \bigoplus_{j=1}^a h_j^*T_{\cX}
\longrightarrow 
\bigoplus_{j=1}^a (T_xX)_{q_j}
\longrightarrow 0 .
\]
The last arrow takes the difference of the two values at each node. It is surjective on global sections, since the bundle on each tail $C_j$ is globally generated and that tail contains just one node. The central summand is a sum of copies of \(\OO_{\PP^1}(-1)\), and the tail summands split into line bundles of nonnegative degree. All these bundles have vanishing \(H^1\). It follows that $$H^1(C,V(-z))=0.$$

A smoothing of this pointed source $C$ is obtained by blowing up \(C_0\times\mathbb A^1\) at the points 
\((q_j,0), j=1,...,a\),
$$\text{Bl}_{\{(q_1,0),...,(q_a,0)\}}(C_0\times \mathbb{A}^1).$$ 
Identify each exceptional line, together with its point on the strict transform of \(C_0\times\{0\}\), with 
\((C_j,p_j)\). The central fiber is then \(C\), both markings extend as constant sections away from the blowup centers, and every other fiber is \(\PP^1\). Consider maps from the fibers of this family to \(\cX\) taking the section through 
\(z\) to \(x\). The stack of such maps is algebraic and locally of finite presentation near \(h\). Indeed, we can use the stack of stable maps from \cite{AV}, with an identification of the marked source with the prescribed fiber, and impose the evaluation at \(z\). For a small extension with ideal \(I\), smoothness of \(\cX\) identifies the obstruction to lifting a map over a prescribed lift of its source with a class in \(H^1(C,V)\otimes I\), see e.g. \cite[Theorem~1.5]{Olsson}. With the value at \(z\) prescribed, the local lifts are chosen to agree there, so their differences take values in \(V(-z)\otimes I\). The obstruction group is thus \(H^1(C,V(-z))\otimes I=0\). The infinitesimal criterion gives smoothness over \(\mathbb A^1\) at \(h\), and the family of curves therefore admits a lift through \(h\) after an \'etale base change near \(0\). See Figure \ref{fig:campana-comb-smoothing} for a depiction.
 
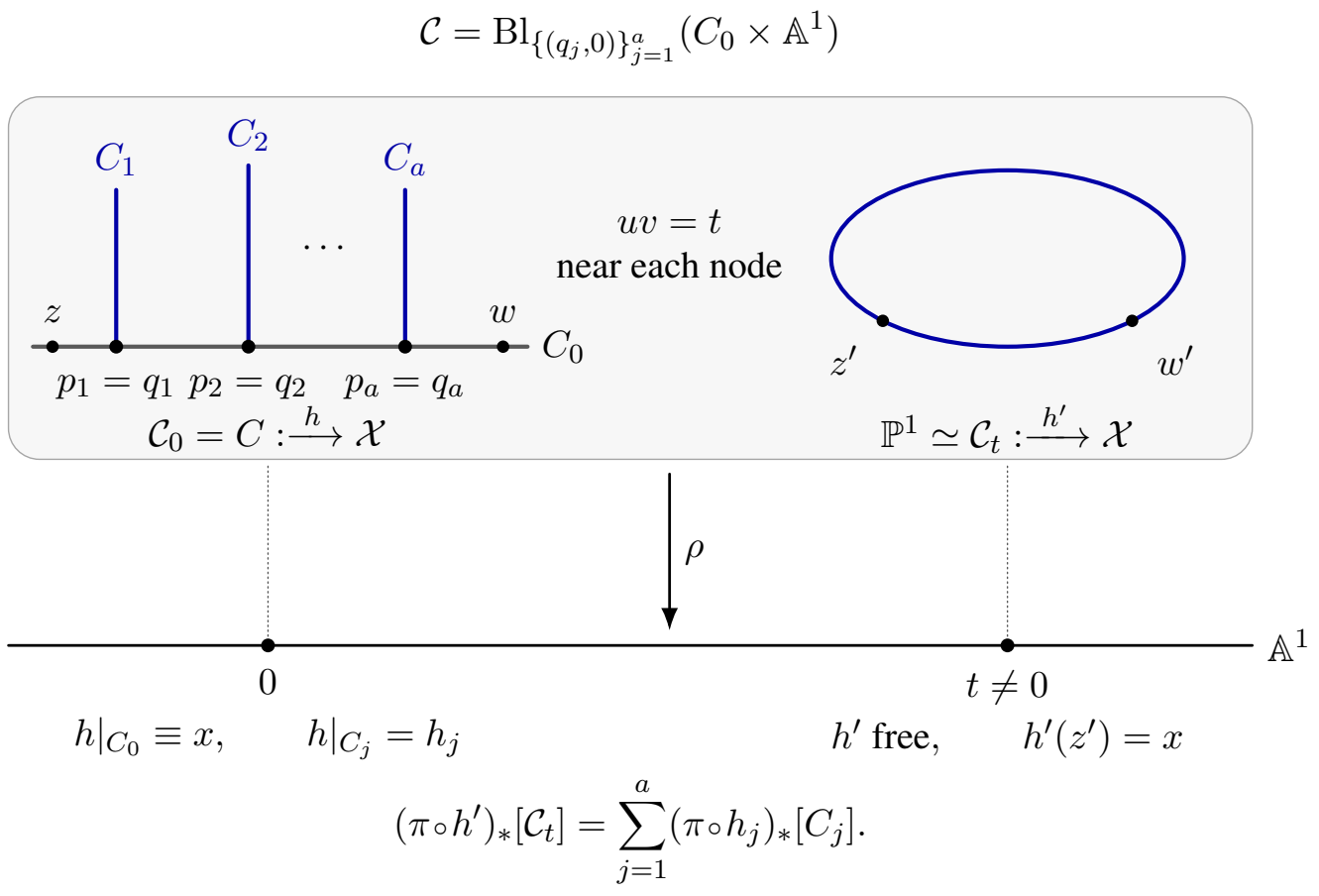
\begin{figure}[htbp]
\centering
\resizebox{0.50\textwidth}{!}{%
\begin{tikzpicture}[x=1cm,y=1cm,font=\small,
  line cap=round,line join=round,>=Latex]

  \filldraw[fill=black!3,draw=black!35,rounded corners=9pt]
    (-0.4,0.65) rectangle (12.3,4.35);
  \node at (5.95,4.95)
    {$\displaystyle\mathcal C=
      \operatorname{Bl}_{\{(q_j,0)\}_{j=1}^{a}}
      (C_0\times\mathbb A^1)$};

  \draw[very thick,black!65] (-0.15,1.8)--(4.9,1.8);
  \node[right] at (4.9,1.8) {$C_0$};

  \foreach \s/\j/\height in
    {0.7/1/3.4,2.05/2/3.65,3.65/a/3.4}{
    \draw[very thick,blue!65!black]
      (\s,1.8)--(\s,\height);
    \node[above,blue!65!black] at (\s,\height) {$C_{\j}$};
    \fill (\s,1.8) circle (2pt);
    \node[below=4pt] at (\s,1.8) {$p_{\j}=q_{\j}$};
  }
  \node at (2.85,2.8) {$\cdots$};
  \fill (0.05,1.8) circle (1.8pt)
    node[above=3pt] {$z$};
  \fill (4.65,1.8) circle (1.8pt)
    node[above=3pt] {$w$};
  \node at (2.25,0.98) {$\mathcal C_0=C: \xrightarrow{\ h\ }\mathcal{X}$};

  \begin{scope}[shift={(9.8,2.7)}]
    \draw[very thick,blue!65!black]
      (0,0) ellipse [x radius=1.8cm,y radius=0.9cm];

    \fill ({1.8*cos(225)},{0.9*sin(225)}) circle (1.8pt)
      node[below left=3pt] {$z'$};
    \fill ({1.8*cos(315)},{0.9*sin(315)}) circle (1.8pt)
      node[below right=3pt] {$w'$};

  \end{scope}
  \node at (9.8,0.98) {$\mathbb P^1\simeq\mathcal C_t: \xrightarrow{\ h'\ }\mathcal{X}$};

  \node[align=center] at (6.35,2.85)
    {$uv=t$\\near each node};

  \draw[thick] (-0.4,-1.25)--(12.3,-1.25)
    node[right] {$\mathbb A^1$};
  \draw[->,thick] (6.35,0.5)--(6.35,-1.1)
    node[midway,right] {$\rho$};

  \draw[densely dotted,black!65]
    (2.25,0.58)--(2.25,-1.25);
  \draw[densely dotted,black!65]
    (9.8,0.58)--(9.8,-1.25);
  \fill (2.25,-1.25) circle (2pt)
    node[below=3pt] {$0$};
  \fill (9.8,-1.25) circle (2pt)
    node[below=3pt] {$t\ne0$};

  \node at (2.25,-2.2)
    {$h|_{C_0}\equiv x,\qquad h|_{C_j}=h_j$};
  \node at (9.8,-2.2)
    {$h'\text{ free},\qquad h'(z')=x$};

  \node at (5.95,-3.15)
    {$\displaystyle
      (\pi\mathbin{\!\scriptstyle\circ\!} h')_*[\mathcal C_t]
      =\sum_{j=1}^{a}(\pi\mathbin{\!\scriptstyle\circ\!} h_j)_*[C_j].$};

\end{tikzpicture}%
}

\caption{\small Two fibers of the family $\rho:\mathcal C\to\mathbb A^1$. The special fiber is the comb $\mathcal C_0=C_0\cup C_1\cup\cdots\cup C_a$; every fiber over $t\ne0$ is a smooth rational curve. The markings $z,w$ extend to disjoint sections, with values $z',w'$ on the displayed smooth fiber. After an \'etale base change near $0$, the map $h$ extends with the image of the first marking fixed at $x$. A nearby smooth fiber carries the free map $h'$. The shaded region represents the total space schematically.}
\label{fig:campana-comb-smoothing}
\end{figure}

Write \(h'\) for the map on a nearby smooth fiber and \(z'\) for its marking. Semicontinuity implies  
\(H^1(\PP^1,(h')^*T_{\cX}(-z'))=0\). On \(\PP^1\), this vanishing says that every rank $1$ summand of \((h')^*T_{\cX}\) has nonnegative degree. Thus \(h'\) is free on the stack $\cX$. For any line bundle on \(X\), its degree on the smooth fiber equals its degree on \(C\), which is the sum of its degrees on the tails. This proves the assertion about the
coarse curve class.
\end{proof}

In each of our finitely many covering families, the free members sweep out a dense constructible subset of \(X\). Choose \(x\in X\setminus D\) in the common dense open contained in these subsets and take one free member of each family through \(x\). Applying Lemma~\ref{lem:smooth-sum}, we obtain a free map 
$$h':\PP^1\to \cX$$ 
of coarse curve class \(\beta\).

\subsection{Boundary contacts and logarithmic freeness}
We now allow the image of the marking to move in order to obtain transversality to the root divisors.

\begin{lemma}\label{lem:free-contacts}
Let \(h':\PP^1\to\cX\) be free on the stack, with nonconstant coarse map meeting \(X\setminus D\). A general deformation \(g\) of \(h'\) avoids the deeper strata of the reduced root divisor \(\mathcal R=\sum_i\mathcal R_i\) and meets its components transversely. After marking all boundary contact points, \(f=\pi\mathbin{\!\scriptstyle\circ\!} g\) is a free Campana curve. Every contact of \(f\) with \(D_i\) has order \(m_i\).
\end{lemma}

\begin{proof}
Choose a point \(z'\in\PP^1\) with \(h'(z')\) outside the root divisor. Take an irreducible \'etale scheme chart \(M\) in the stack of parametrized maps from \(\PP^1\) near \(h'\). Shrink it so that every map \(g\) in this chart is free and \(g(z')\) is outside the boundary divisor. This chart $M$ is smooth because \(H^1(\PP^1,g^*T_{\cX})=0\). Moreover, the universal evaluation $$e:\PP^1\times M\to\cX$$
is smooth. Indeed, if the value at a point \(\mathsf p\) is prescribed, the obstruction group becomes
\(H^1(\PP^1,g^*T_{\cX}(-\mathsf p))\), as in the proof of Lemma~\ref{lem:smooth-sum}. This group vanishes for every \(\mathsf p\) by freeness. The same lifting argument allows \(\mathsf p\) to vary as a marked section.

For \(i\ne j\), the inverse image of \(\mathcal R_i\cap\mathcal R_j\) under \(e\) is empty or has codimension $2$. Its image under the proper projection \(\PP^1\times M\to M\) is closed of dimension at most \(\dim M-1\). Remove these images and put $$W_i=e^{-1}(\mathcal R_i).$$
Each nonempty \(W_i\) is a smooth divisor of dimension \(\dim M\). Every fiber of \(W_i\to M\) is finite, since the curve contains a point outside the boundary divisor. The morphism $W_i\to M$ is proper, hence finite, and each of its irreducible components dominates \(M\) by dimension count. In characteristic zero, it is generically \'etale on each component. Removing the images of the non-\'etale loci makes the \(W_i\) pairwise disjoint finite \'etale covers. In particular, their fibers are reduced, which is precisely transversality to the root divisors.

Write $$F=\pi\mathbin{\!\scriptstyle\circ\!} e: \PP^1\times M\to X.$$ 
The equality of relative Cartier divisors \(F^*D_i=m_iW_i\) shows that each contact with \(D_i\) has order \(m_i\). If the coarse curve class is \(\beta\), the degree of \(W_i\to M\) is \(d_i=(D_i\cdot\beta)/m_i\). We may label the points of every \(W_i\) after taking the product over \(M\) of the finite \'etale covers parametrizing orderings of its fibers. On this base change, \(W_i\)'s are unions of disjoint marked sections. Equip the domain of (the base-changed) $F$ with their divisorial log structure and the base with the trivial log structure. A local equation for \(D_i\) now pulls back to a unit times the \(m_i\)-th power of a local equation of the marked section at each of its contacts. Thus \(F\) defines a family of nondegenerate log maps with fixed contact data. Each coarse map is nonconstant, so the marked maps are stable.

It remains to check freeness in the sense of \cite[Definition~4.2]{CLT}. Fix a member \(g\) and let \(S\) be the
reduced divisor of contact markings on \(P=\PP^1\). Under \eqref{eq:kummerlog}, the logarithmic normal complex of 
\(f=\pi\mathbin{\!\scriptstyle\circ\!} g\) is the top row of the following diagram:
\begin{equation*}
\xymatrix{
T_P(-S)\ar[d]\ar[r]^(0.4){d g^{\log}}& g^*\!\left(T_{\cX}(-\log\mathcal R)\right)\ar[d]\\
T_P\ar[r]^{d_g} & g^*T_{\cX}.
}    
\end{equation*}
Both rows are placed in degrees \(-1,0\). Away from \(S\), the vertical arrows are isomorphisms. At a contact \(\mathsf p\) with \(\mathcal R_i\), choose a source parameter \(t\) and a normal coordinate \(u\) for \(\mathcal R_i\) in an \'etale chart. Transversality implies \(u\mathbin{\!\scriptstyle\circ\!} g=t v(t)\), where \(v(0)\ne0\). The vertical arrows are injective, given in local frames by multiplication by \(t\) on the left and by \(\operatorname{diag}(t v(t),1,\ldots,1)\) on the right. The quotient at \(\mathsf p\) is the two-term complex
\[
 T_{\mathsf{p}}P\ \xrightarrow{\ d g^{\perp}_{\mathsf{p}}\ }\ (g|_{\mathsf{p}})^*N_{\mathcal R_i/\cX}.
\]
Its differential is multiplication by \(v(0)\), hence is an isomorphism. The two rows are therefore quasi-isomorphic on 
\(P\).

The coarse map is nonconstant and meets \(X\setminus D\), so in characteristic zero \(dg\) is generically nonzero. Its source is a line bundle, hence \(dg\) is injective as a map of sheaves. Both normal complexes are thus represented by the sheaf \(N_g=\operatorname{coker}(dg)\), including any torsion at ramification points. Tensoring the defining
sequence of \(N_g\) with \(\OO_P(-1)\), we obtain
$$0\to T_P(-1)\xrightarrow{\ d_g\ }g^*T_{\cX}(-1)\to N_g(-1)\to 0.$$ 
Since \(H^1(P,g^*T_{\cX}(-1))=0\) and \(H^2(P,T_P(-1))=0\), the associated long exact sequence gives 
\(H^1(P,N_g(-1))=0\). The quasi-isomorphism identifies this group with the corresponding hypercohomology of the logarithmic normal complex. This proves the desired freeness.
\end{proof}

\subsection{Conclusion of the proof of Theorem~\ref{thm:main}}
We apply Lemma~\ref{lem:free-contacts} to the free map of class \(\beta\) constructed above. The coarse curve class of a deformation is again \(\beta\). Its intersection with each \(D_i\) is positive because \(\beta\) is interior and the class \([D_i]\) is nonzero and effective. Thus \(f\) meets every \(D_i\), with the exact orders as claimed. By \cite[Definition~5.9]{CLT}, this free curve of interior nef class makes \((X,\Delta)\) strongly Campana uniruled.

\cite[Theorem~1]{Zhang} applies to the klt pair \((X,\Delta)\) and gives rational connectedness of \(X\). Since the ground field is \(\mathbb C\), \cite[Corollary~6.7]{CLT} applies. Its conclusion, in the terminology of \cite[Definitions~4.1(2) and~5.6(2)]{CLT}, is a family with fixed Campana contact data and dominant two-point evaluation.

\subsection{Surface examples}\label{sec:example}
We present here two concrete examples to illustrate Theorem \ref{thm:main}.

\subsubsection{A weak log Fano pair on a Hirzebruch surface}
Let 
$$q:X=\mathbb F_2=\PP_{\PP^1}(\OO\oplus\OO(2))\to\PP^1$$ 
be the projective bundle of $1$-dimensional quotients. Write \(F\) for the class of a fiber and \(H=c_1(\OO_X(1))\). The quotients onto \(\OO\) and \(\OO(2)\) give disjoint sections \(S\) and \(D\), respectively. The map from the trivial summand to the tautological quotient vanishes precisely along \(D\), with multiplicity one on each fiber. Applying the same description to the other summand gives
\[
D\sim H,\qquad S\sim H-2F.
\]
In particular \(H\cdot S=0\), \(H\cdot F=1\), and \(H^2=2\). The bundle \(\OO\oplus\OO(2)\) is globally generated, hence so is \(\OO_X(1)\). The projective bundle formula for the canonical divisor gives 
$$K_X\sim-2H+q^*(K_{\PP^1}+c_1(\OO(2)))\sim-2H.$$
For \(m\ge2\), it follows that
\[
-\bigl(K_X+(1-1/m)D\bigr)\sim_{\mathbb Q}(1+1/m)H .
\]
This divisor is nef and has positive square, so it is big. Its degree on \(S\) is $0$. Thus the weak log Fano hypothesis of Theorem~\ref{thm:main} applies to this pair.

\begin{figure}[htbp]
\centering
\resizebox{0.50\textwidth}{!}{%
\begin{tikzpicture}[
  x=0.95cm,y=0.90cm,font=\small,
  line cap=round,line join=round,>=Latex]

  \filldraw[fill=black!3,draw=black!30,rounded corners=5pt]
    (0,0) rectangle (5.4,3.1);
  \node at (2.7,3.55) {$X=\mathbb F_2$};
  \draw[thick,blue!65!black] (0.15,2.4)--(5.25,2.4);
  \node[above right,blue!65!black] at (0.15,2.4)
    {$D\sim H,\quad D^2=2$};
  \draw[thick] (0.15,0.7)--(5.25,0.7);
  \node[below right] at (0.15,0.7)
    {$S\sim H-2F,\quad S^2=-2$};
  \draw[thick,black!55] (4.7,0.05)--(4.7,3.05);
  \node[right] at (4.7,1.55) {$F$};
  \fill (4.7,0.7) circle (1.6pt);
  \fill[blue!65!black] (4.7,2.4) circle (1.6pt);
  \draw[->] (4.7,-0.12)--node[right] {$q$} (4.7,-1.12);
  \draw[thick] (0,-1.3)--(5.4,-1.3);
  \fill (4.7,-1.3) circle (1.6pt);
  \node[right] at (5.4,-1.3) {$\mathbb P^1$};

  \begin{scope}[xshift=7.2cm]
    \node at (1.6,3.55) {Local contact model};
    \node at (1.6,2.9) {$C:\ z=t^m$};
    \draw[thick,blue!65!black] (-0.15,1.1)--(3.45,1.1);
    \node[below right,blue!65!black] at (2.15,1.1)
      {$D:\ z=0$};
    \draw[very thick,red!65!black]
      plot[domain=-1.1:1.1,samples=60]
      ({1.6+1.2*\x},{1.1+0.95*\x^3});
    \fill (1.6,1.1) circle (1.8pt);
    \node[above left] at (1.6,1.1) {$p$};
    \node at (1.6,-0.65) {$i_p(C,D)=m$};
    \node at (1.6,-1.3) {one of the $d$ contacts};
  \end{scope}

  \node at (5.6,-2.15)
    {$C\equiv S+aF,\quad a=md>2,\quad
      C\cdot S=a-2,\quad C\cdot F=1.$};
\end{tikzpicture}%
}
\caption{\small The Hirzebruch surface example. The sections $S$ and $D$ are disjoint, and each meets a fiber once. The constructed section $C=f(\mathbb P^1)$ has $d$ contacts of order $m$ with $D$. The inset gives the local equation at each contact; its sketch uses $m=3$. The global curve $C$ is not drawn in the left panel.}
\label{fig:campana-hirzebruch-surface}
\end{figure}
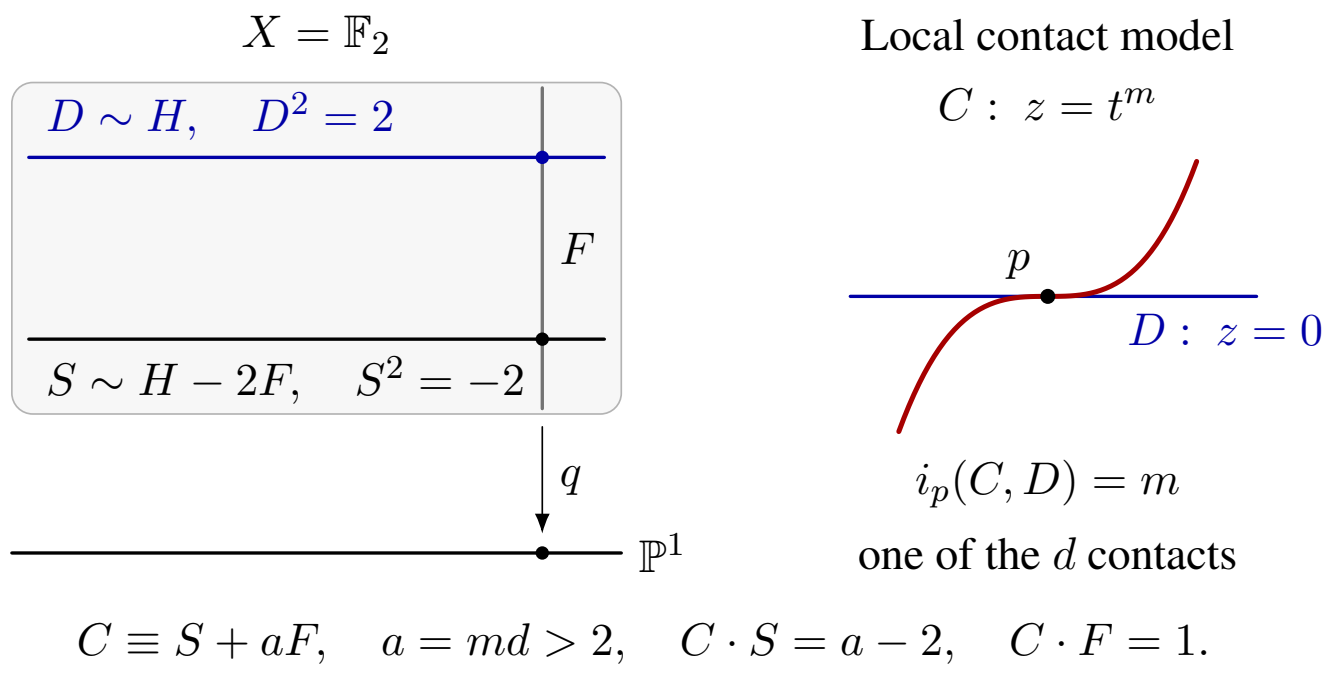

We construct a section of the ruling whose lift to the root stack \(\pi:\cX\to X\) is very free. Choose \(d\ge1\) such that \(a=md>2\), and choose a section \(u\in H^0(\PP^1,\OO(d))\) with \(d\) distinct zeros. Since \(\OO(a-2)\) is globally generated, vanishing at any one of these points imposes a proper linear condition on its space of sections. We may therefore take \(v\in H^0(\PP^1,\OO(a-2))\) which is nonzero at every zero of \(u\). The homomorphism
\[
(u^m,v):\OO\oplus\OO(2)\longrightarrow\OO(a)
\]
is surjective at each point and defines a section \(f:\PP^1\to X\). The equation of \(D\) pulls back to \(u^m\), hence \(f^*D=m\operatorname{div}(u)\). The pair \((\OO(d),u)\) gives a lift \(g:\PP^1\to\cX\) with simple contacts with the reduced root divisor.

To identify the numerical class of \(f\), note that \(S,F\) form a basis of \(N^1(X)_{\mathbb R}\). For any irreducible curve \(C\) on this surface, we have
\[
C\equiv(C\cdot F)S+(C\cdot H)F.
\]
Both coefficients are nonnegative since \(F\) and \(H\) are nef. As \(S\) and \(F\) are themselves effective, they generate the pseudoeffective divisor cone. Now \(f^*\OO_X(1)=\OO(a)\) and 
\(q\mathbin{\!\scriptstyle\circ\!} f=\mathrm{id}\), so \(f_*[\PP^1]=S+aF\). Its intersections with 
\(S\) and \(F\) are \(a-2\) and \(1\) respectively. This proves that the class is interior nef.

The equality \(q\mathbin{\!\scriptstyle\circ\!} \pi\mathbin{\!\scriptstyle\circ\!} g=\mathrm{id}_{\PP^1}\) splits the differential \(g^*T_{\cX}\to T_{\PP^1}\). Its kernel is a line bundle, and \eqref{eq:rootcanonical} gives \(\deg g^*T_{\cX}=(1+1/m)a=(m+1)d\). Consequently
\[
g^*T_{\cX}\simeq\OO(2)\oplus\OO((m+1)d-2).
\]
Both summands on the right have positive degree. Thus \(g\) is very free on the stack. Fix distinct \(t_1,t_2\in\PP^1\) away from the contacts. The above splitting gives \(H^1(\PP^1,g^*T_{\cX}(-t_1-t_2))=0\), so evaluation at these two points is smooth at \(g\) by the pointed deformation argument of Lemma~\ref{lem:smooth-sum}. Both evaluation values lie in \(X\setminus D\). The contacts stay simple upstairs on an open neighborhood of \(g\) and can be labeled after finite \'etale base change as in Lemma~\ref{lem:free-contacts}. Thus the coarse family has dominant two-point evaluation and all contacts have order \(m\).

\subsubsection{A conic and two lines in the plane}
Let \(B\subset\PP^2\) be a smooth conic and \(L_1,L_2\subset \PP^2\) distinct lines such that \(D=B+L_1+L_2\) is SNC. Give \(B,L_1,L_2\) the multiplicities \(a,b,c\ge2\), respectively. If \(H\) is the class of a line, then 
\(K_{\PP^2}=-3H\) and \(B\sim 2H\). Hence
\[
-(K_{\PP^2}+\Delta)\sim_{\Q}\lambda H,\quad \text{where}\quad \lambda=\frac{2}{a}+\frac1b+\frac1c-1 .
\]
The weak log Fano condition for $(\PP^2, \Delta)$ is equivalent to \(\lambda>0\). In this case, Lemma~\ref{lem:covering} applied to the root stack of the pair, already produces a free curve of interior nef class. Indeed, the class
of its nonconstant coarse map has positive intersection with \(H\), and \(H\) generates the pseudoeffective divisor cone.

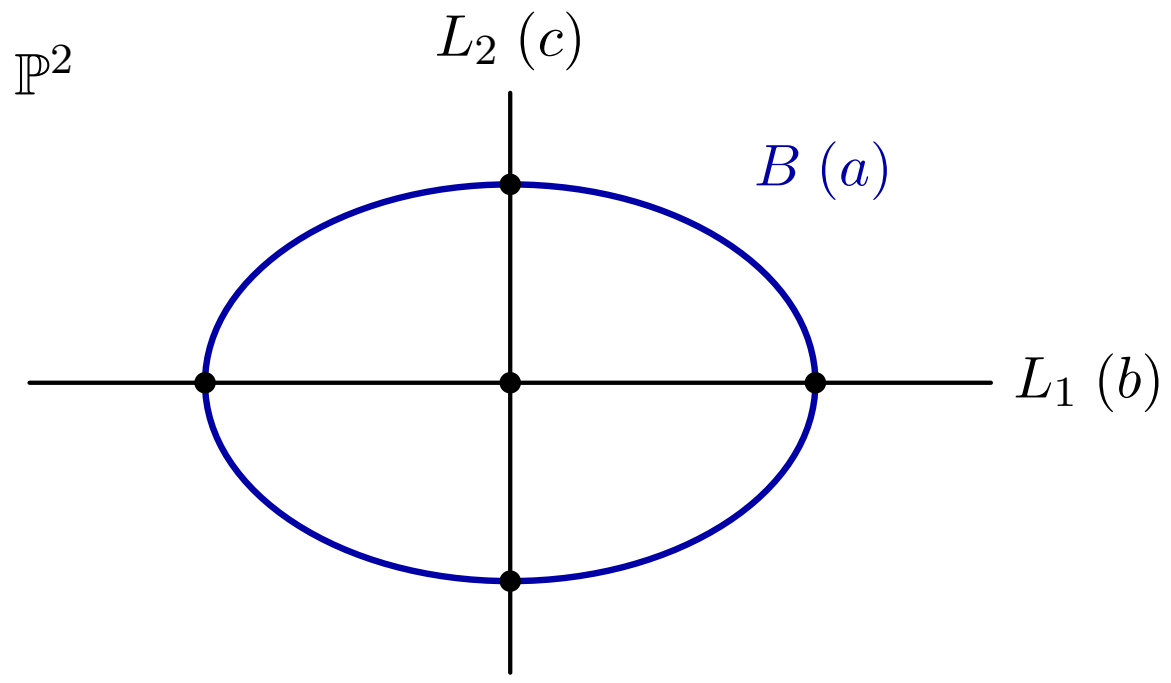
\begin{figure}[htbp]
\centering
\resizebox{0.50\textwidth}{!}{%
\begin{tikzpicture}[
  x=1cm,y=1cm,font=\small,
  line cap=round,line join=round]

  \node at (-3.05,2.05) {$\mathbb P^2$};
  \draw[very thick,blue!65!black]
    (0,0) ellipse [x radius=2cm,y radius=1.3cm];
  \draw[thick] (-3.15,0)--(3.15,0);
  \draw[thick] (0,-1.9)--(0,1.9);
  \node[right] at (3.15,0) {$L_1\ (b)$};
  \node[above] at (0,1.9) {$L_2\ (c)$};
  \node[above right,blue!65!black] at (1.45,1.05)
    {$B\ (a)$};

  \foreach \p in {(-2,0),(2,0),(0,-1.3),(0,1.3),(0,0)}
    \fill \p circle (2pt);

  \node at (0,-2.5)
    {$\displaystyle
      \Delta=\left(1-\frac1a\right)B
      +\left(1-\frac1b\right)L_1
      +\left(1-\frac1c\right)L_2$};
  \node at (0,-3.5)
    {$\displaystyle
      -(K_{\mathbb P^2}+\Delta)\sim_{\mathbb Q}\lambda H,
      \qquad \lambda=\frac2a+\frac1b+\frac1c-1>0.$};
\end{tikzpicture}%
}
\caption{\small A smooth conic and two lines in the plane. The numbers in parentheses are the orbifold multiplicities. The five marked points are transverse double intersections; there are no triple intersections.}
\label{fig:campana-conic-two-lines}
\end{figure}

Let \(g:\PP^1\to\cX\) be any free map to this root stack whose coarse map \(f\) has degree \(d>0\) and meets the boundary divisor complement. The divisors \(f^*B,f^*L_1,f^*L_2\) have degrees \(2d,d,d\). The root conditions therefore give \(a\mid 2d\), \(b\mid d\), and \(c\mid d\). By \eqref{eq:rootcanonical}, the degree of \(g^*T_{\cX}\) is 
\(\lambda d\). Its summands have nonnegative degrees, and the nonzero differential \(dg:\OO(2)\to g^*T_{\cX}\) forces one degree to be at least $2$. Thus \(d\) is divisible by \(\operatorname{lcm}(a/\gcd(a,2),b,c)\) and satisfies 
\(d\ge 2/\lambda\). The resulting necessary conditions for two triples are recorded below.
\[
\begin{array}{c|c|c|c}
 (a,b,c)&\lambda&\text{divisibility of }d&\text{lower bound for }d\\ \hline
 (3,4,11)&1/132&132\mid d&d\ge264\\
 (3,5,7)&1/105&105\mid d&d\ge210
\end{array}
\]
Note that this example is covered by \cite[Corollary 1.8]{CT}.


\part{Klt case}\label{part:klt}
Here we treat the klt case of Theorem \ref{thm:main_gen}.

\section{Overview: klt case}\label{cklt:appendix}

Let $X$, $D_1,\ldots,D_s$, $D=\sum_i D_i$, and $\Delta=\sum_i(1-1/m_i)D_i$ be as in Section \ref{sec:result}. We allow $s=0$. On the log-smooth locus of $(X,D)$ we use the contact orders and the two-point evaluation of Section~\ref{sec:result}.

Denote by $X_{\text{reg}}\subset X$ the smooth locus. Define the {\em basic open} 
\begin{equation}\label{eqn:basic_open}
X^\circ\subset X_{\text{reg}} 
\end{equation}
to be the complement in $X_{\mathrm{reg}}$ of the singular loci of $D_1,...,D_s$ and pairwise intersections $D_i\cap D_j, 1\leq i\neq j\leq s$. Let 
\begin{equation}\label{eqn:root_stack_over_basic_open}
\pi:\mathcal X\to X^\circ
\end{equation}
be the stack of $m_i$-th roots along $D_i$, $i=1,...,s$. Recall that a map from $\mathbb P^1$ to $\mathcal X$ with nonconstant coarse image is {\em free} if the pullback of $T\mathcal{X}$ is a globally generated vector bundle on $\PP^1$ and {\em very free} if that bundle is ample. As before, $\operatorname{Nef}_1(X)$ is the dual of the cone of pseudoeffective numerical Cartier divisor classes.

\begin{theorem}[Campana rational connectedness]\label{cklt:thm:dlt-descent}
Suppose $(X,\Delta)$ is klt and $-(K_X+\Delta)$ is nef and big. The root stack $\pi:\mathcal{X}\to X^\circ$ in (\ref{eqn:root_stack_over_basic_open}) admits a very free map $\mathbb{P}^1\to \mathcal{X}$ whose coarse class lies in $\operatorname{int}\operatorname{Nef}_1(X)$. For a general deformation of this very free map, the coarse map $f$ meets each $D_i$, and
\[
\operatorname{ord}_p f^*D_i=m_i \qquad\text{for every }p\in f^{-1}(D_i).
\]
With its contact points marked, $f$ is a very free logarithmic curve for this contact data. The coarse maps form a family with fixed contact data and dominant two-point evaluation. Thus, $(X,\Delta)$ is Campana rationally connected.
\end{theorem}

Here logarithmic very freeness of $f$ means that the logarithmic normal complex $N_f^{\log}$ satisfies $\mathbb H^1(\mathbb{P}^1,N_f^{\log}(-2))=0$, see \cite[Definition~4.2]{CLT}.

We spell out the $D=\emptyset$ case.
\begin{corollary}[Very free curves in the smooth locus]\label{cklt:cor:weak-fano-smooth}
Let $X$ be a positive-dimensional normal projective klt variety with $-K_X$ nef and big. Then there is a nonconstant morphism $f:\mathbb{P}^1\to X_{\mathrm{reg}}$ with $f^*T_{X_{\mathrm{reg}}}$ ample. Its deformations in $X_{\mathrm{reg}}$ have dominant two-point evaluation.
\end{corollary}

Index-one components of a specified reduced divisor may also be included. They impose contact order at least one. A general deformation avoids their singular loci and pairwise intersections with the other components, and is transverse at its contacts, by smoothness of universal evaluation. After a finite \'etale base change these additional contacts may be marked as well.

\subsection{Proof strategy}
Choose a log resolution $p:V\to X$ and take roots along the strict boundary and the exceptional divisors. The weak log Fano assumption for $(X,\Delta)$ need not hold for the pair constructed from taking roots on $V$. In Section~\ref{cklt:sec:divisor}, we handle this by applying Lemma \ref{lem:arrangement} on a smooth model of the base. 

Our basic strategy is to establish a modified version of Lemma~\ref{lem:covering} which applies to the pair $(X,\Delta)$. As the klt case is more delicate, we cannot  apply Lemma~\ref{lem:covering} directly. To circumvent this difficulty, we begin  with an auxiliary construction   involving a projective bundle, a  resolution of that bundle, and root stacks. This construction, together with a careful  discrepancy comparison, allows us to use Lemma~\ref{lem:covering} to construct free curves in relation to divisors on $V$. The upshot of Section~\ref{cklt:sec:divisor} is that the free quotient is a point and that we can choose a very free curve whose class projects to the interior of $\operatorname{Nef}_1(X)$.
See Proposition~\ref{cklt:prop:singular-divisor-test} and  Proposition~\ref{cklt:prop:initial} for more details.

The next issue we need to address  is that the curve can still meet an exceptional divisor of $p$. 
A priori, the free curve on $V$ may meet the exceptional divisors and its pushforward to $X$ may fail to  retain  the desired properties.  
To deal with this, we run a relative minimal model program. 
We can treat the exceptional divisors one at a time. 
At a divisorial step, adjunction yields a lower-dimensional log Fano pair on a general fiber of that divisor. By induction, we obtain suitable curves in the fiber.  Our idea is to attach  fiber curves  as teeth  to the comb. The intersection-theoretic consideration provides a candidate of a curve class avoiding the exceptional divisor. We then smooth the whole curve. More precisely, we choose the root order so that a cover of a fiber curve has normal degree $-1$, and attach one such curve at each transverse contact. Each tooth cancels one unit of intersection. Moving the attachments allows us to smooth the nodes with two evaluations prescribed. 
The smoothed curve misses the exceptional divisor and has the same class on $X$. 
The details of canceling transverse contact points 
and finding suitable teeth via orbifold adjunction are in Section~\ref{cklt:sec:divisorialcontraction}. 
We conclude the proof in Section \ref{cklt:sec:descent}.

\section{Constructing free curves relative to a divisor}\label{cklt:sec:divisor}
We work in this situation of Theorem~\ref{cklt:thm:dlt-descent}. Set $A=-(K_X+\Delta)$ and fix a projective log resolution 
\begin{equation}\label{eqn:log_resolution}
p:V\to X.
\end{equation}
Write 
$$\nu:\mathcal V\to V$$ 
for the stack of $m_i$-th roots $m_i$ along the strict transforms of $D_i$ and $l_k\ge 2$-th roots  along the exceptional prime divisors $E_k\subset V$ of $p$. Here $l_k\geq 2$ are arbitrarily chosen integers.

\begin{proposition}[Curves against a divisor]\label{cklt:prop:singular-divisor-test}
Let $L$ be an effective divisor on $V$ such that $E=p_*L$ is nonzero and $\mathbb Q$-Cartier. Then there is a dominating family of free maps $f:\mathbb{P}^1\to\mathcal V$ with nonconstant coarse map and $\deg(\nu\mathbin{\!\scriptstyle\circ\!} f)^*L>0$. 
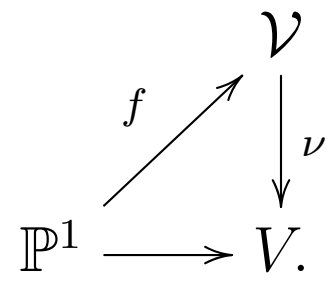
\begin{equation*}
\xymatrix{
 &\mathcal{V}\ar[d]^\nu\\
\mathbb{P}^1\ar[r]\ar[ur]^f& V.
}    
\end{equation*}

For $L=p^*E$ with $E$ a nonzero effective Cartier divisor on $X$, the auxiliary construction below gives a projective bundle $q:Y\to X$, a resolution $\mu:W\to Y$, and roots $s:\mathcal W\to W$ such that
\[
 K_{\mathcal W}\sim_{\mathbb Q}s^*(-\varepsilon(q\mathbin{\!\scriptstyle\circ\!}\mu)^*E+R),\qquad R\ge0\ \text{and }\mu\text{-exceptional}.
\]
The general curve $g$ in the covering family on $\mathcal W$ satisfies $\varepsilon\deg(q\mathbin{\!\scriptstyle\circ\!}\mu\mathbin{\!\scriptstyle\circ\!} s\mathbin{\!\scriptstyle\circ\!} g)^*E\ge2$.
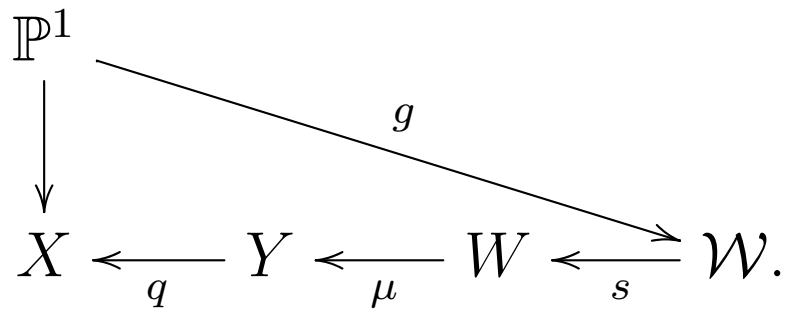
\begin{equation*}
\xymatrix{
\PP^1\ar[d]\ar[drrr]^g& & & \\
X & Y\ar[l]^q& W\ar[l]^{\mu}&\mathcal{W}\ar[l]^{s}.
}    
\end{equation*}

\end{proposition}

\begin{proof}
We adapt the proof of Proposition~\ref{prop:testing}. The divisor $L$ is on $V$, but the weak log Fano is on $X$. We first construct the auxiliary pair over $X$ and then compare the two divisors on a common model.

\smallskip
\noindent\emph{The auxiliary pair.}
Write $$A\sim_{\mathbb Q}H+G_0$$ 
with $H$ ample and $G_0$ effective and $\mathbb Q$-Cartier. Following the proof of Lemma~\ref{lem:perturbation}, fix a log resolution $h:Z\to X$ of 
$\operatorname{Supp}(\Delta+E+G_0)$ and write
\[
K_Z+\Gamma_\varepsilon=h^*(K_X+\Delta+\varepsilon E).
\]
Let $\Gamma_\varepsilon^+$ denote the coefficientwise positive part of $\Gamma_\varepsilon$. At $\varepsilon=0$, all coefficients of $\Gamma_0^+$ are less than $1$. Choose small positive rational numbers $t$, then $\varepsilon$, so that $\Gamma_\varepsilon^++t h^*G_0$ has coefficients less than $1$ on the fixed SNC support. Setting $T=tG_0$, we see that $(X,\Delta+T)$ is klt and 
$$A-T\sim_{\mathbb Q}(1-t)A+tH$$ 
is ample. We may also require that 
$$P=A-T-\varepsilon E$$ 
is ample. For sufficiently divisible $m\ge3$, $mT$ is Cartier and $mP$ is very ample. Take a general $G_1\in|mP|$ and put $G=mT+G_1$. Then $$G\sim_{\mathbb Q}m(A-\varepsilon E).$$
Applying Bertini to the pullback linear system, we see that $$\left(Z,\Gamma_\varepsilon^++\frac{1}{m}h^*G\right)$$ is klt.

Put $r=m-2$, $w=1-1/m$, and let 
\begin{equation}\label{eqn:proj_bundle_klt}
q:Y=\mathbb{P}_X(\mathcal O_X^{\oplus r}\oplus\mathcal O_X(G))\to X
\end{equation}
be the projective bundle that parametrizes $1$-dimensional quotients. Let $\mathcal H$ be the sum of the hyperplane at infinity and the $r+1$ hyperplanes given by $e_1,\ldots,e_r,\sum_j e_j+s_G$ in the tautological quotient, as in Lemma~\ref{lem:arrangement}. Here $s_G$ is the divisor section. Put $$B=q^*\Delta+w\mathcal H,$$ 
using flat pullback for Weil divisors. With $\xi=c_1(\mathcal O_Y(1))$, we have $$\mathcal H\sim m\xi-q^*G \text{ and } K_{Y/X}\sim-(m-1)\xi+q^*G.$$ 
Thus
\[
 K_Y+B\sim_{\mathbb Q}q^*(K_X+\Delta+G/m)\sim_{\mathbb Q}-\varepsilon q^*E.
\]

Applying Lemma~\ref{lem:arrangement} over the smooth base $Z$ with boundary $\Gamma_\varepsilon^+$ and divisor $h^*G$. On $$Y_Z=Y\times_X Z,$$ 
with projections $q_Z:Y_Z\to Z$ and $h_Y:Y_Z\to Y$, it gives a klt pair with boundary divisor $$q_Z^*\Gamma_\varepsilon^++w h_Y^*\mathcal H.$$ 
Lowering $\Gamma_\varepsilon^+$ to $\Gamma_\varepsilon$ preserves the strict discrepancy inequalities. The relative canonical bundle and all the indicated hyperplanes commute with base change. With compatible canonical divisors this gives
\[
 K_{Y_Z}+q_Z^*\Gamma_\varepsilon+w h_Y^*\mathcal H =h_Y^*(K_Y+B+\varepsilon q^*E).
\]
Equality on a common resolution now shows that $$(Y,B+\varepsilon q^*E)$$ 
is klt.

\smallskip\noindent\emph{The divisor on the resolution.}
Choose a log resolution $\mu:W\to Y$ of $(Y,B)$ dominating $Y\times_X V$. Write $b:W\to V$ for its projection, so $q\mathbin{\!\scriptstyle\circ\!}\mu=p\mathbin{\!\scriptstyle\circ\!} b$. 
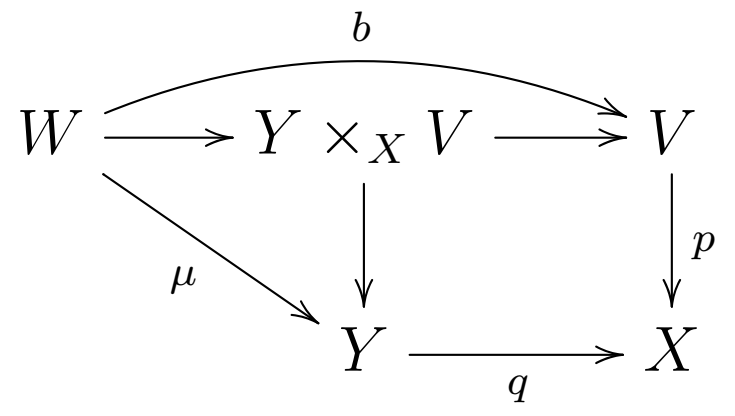
\begin{equation*}
\xymatrix{
W\ar@/^1.20pc/[rr]^{b}\ar[dr]_{\mu}\ar[r] &Y\times_X V\ar[r]\ar[d] &V\ar[d]^p\\
&Y\ar[r]_q& X
}    
\end{equation*}
Let $F_j$ be the $\mu$-exceptional prime divisors on $W$ and $B'$ the strict transform of $B$. Define $a_j,c_j$ by
\[
 K_W+B'=\mu^*(K_Y+B)+\sum_j a_jF_j,\qquad
 (q\mathbin{\!\scriptstyle\circ\!}\mu)^*E-b^*L=\sum_j c_jF_j .
\]
Indeed, $p^*E-L$ is supported over a subset of codimension at least $2$ in $X$, whose inverse image under $q$ still has codimension at least $2$. Put 
$$u_j=\operatorname{ord}_{F_j}((q\mathbin{\!\scriptstyle\circ\!}\mu)^*E).$$ 
Since $b^*L$ is effective, $c_j\le u_j$, and the klt assertion gives $1+a_j-\varepsilon u_j>0$. Choose an integer $N$ divisible by every $m_i,l_k$ and large enough, so that $$1+a_j-1/N-\varepsilon u_j\ge0$$ 
for all $j$.

Let $s:\mathcal W\to W$ be the root stack with the orders of $B'$ and order $N$ along the $F_j$. Its boundary divisor has SNC support. The canonical divisor formula gives
\[
 K_{\mathcal W}\sim_{\mathbb Q}s^*M,\qquad
 M=-\varepsilon b^*L+R,\qquad
 R=\sum_j(1+a_j-1/N-\varepsilon c_j)F_j\ge0 .
\]
The map $b\mathbin{\!\scriptstyle\circ\!} s$ lifts to $\widetilde b:\mathcal W\to\mathcal V$: 
\begin{equation*}
\xymatrix{
\mathcal{W}\ar[d]_s\ar[r]^{\widetilde{b}}\ar[dr] &\mathcal{V}\ar[d]^{\nu}\\
W\ar[r]_b &V. 
}    
\end{equation*}
In fact, the pullback of each strict boundary component on $V$ has a unique non-$\mu$-exceptional component, of multiplicity $1$, rooted with its original order, and its other components are $\mu$-exceptional. The pullback of each $p$-exceptional prime divisor is entirely $\mu$-exceptional. Divisibility of $N$ thus gives the required roots of these Cartier divisors with their sections.

For $J$ very ample on $Y$ and $d=\dim Y$, the projection formula gives 
$$M\cdot(\mu^*J)^{d-1}=-\varepsilon q^*E\cdot J^{d-1}<0.$$ 
A small ample perturbation of $\mu^*J$, followed by rescaling, gives a very ample complete-intersection class of negative $M$-degree. Applying Lemma~\ref{lem:covering} to $\mathcal W$, we obtain a dominating family of free maps $g:\mathbb{P}^1\to\mathcal W$ with $-\deg g^*K_{\mathcal W}\ge2$. A general member is not contained in the support of $R$, so
\[
 \varepsilon\deg(b\mathbin{\!\scriptstyle\circ\!}s\mathbin{\!\scriptstyle\circ\!}g)^*L=-\deg g^*K_{\mathcal W}+\deg(s\mathbin{\!\scriptstyle\circ\!}g)^*R\ge2.
\]
The maps $\widetilde b \mathbin{\!\scriptstyle\circ\!} g$ therefore have nonconstant coarse images and dominate $\mathcal V$ as a family. After reparametrization, evaluation at $0$ is dominant. By generic smoothness, sections of the general pullback tangent bundle span its fiber at $0$, which excludes negative summands and proves freeness. When $L=p^*E$, every $c_j$ vanishes and $b^*L=(q\mathbin{\!\scriptstyle\circ\!}\mu)^*E$, giving the additional assertions.
\end{proof}

\begin{proposition}[The initial very free curve]\label{cklt:prop:initial}
If $X$ is $\mathbb Q$-factorial, $\mathcal V$ admits a very free map whose projected class lies in the interior of $\operatorname{Nef}_1(X)$.
\end{proposition}

\begin{proof}
Let 
\begin{equation}\label{eqn:mfr_quot}
\mathcal V\dashrightarrow Q
\end{equation}
be the maximal free rational quotient \cite[Theorem~1.3, Corollary~2.4(ii) and Proposition~2.6]{cklt:Starr}. (\ref{eqn:mfr_quot}) exists since $\mathcal V$ is smooth and proper, and Proposition~\ref{cklt:prop:singular-divisor-test} gives free curves by using the pullback of an effective ample divisor. A general member of every free component is contracted wherever the quotient is defined, and one free component joins two general points in a general fiber.

If $\dim Q>0$, take a nonconstant rational function on $Q$ and pull it back to $h\in\mathbb{C}(V)=\mathbb{C}(X)$. Write $L$ and $P$ for its zero and pole divisors on $V$. The divisor $p_*L$ on $X$ is nonzero, since otherwise $1/h$ would be regular on the normal projective variety $X$; $p_*L$ is $\mathbb Q$-Cartier because $X$ is $\mathbb Q$-factorial. Applying Proposition~\ref{cklt:prop:singular-divisor-test} to this $L$, we see that there is a free component of positive $L$-degree. Its general member misses $|L|\cap|P|$ on the coarse space. Indeed, free universal evaluation is smooth \cite[Lemma~2.1]{cklt:Starr}, so the inverse image of a codimension-$2$ locus in the universal curve cannot dominate the parameter space. Consequently, along this member, $h$ defines a morphism to $\mathbb{P}^1$. It is constant by the quotient property, with general value outside $\{0,\infty\}$ by dominance. This contradicts positivity of $L$-degree. Hence $Q$ is a point.

The joining component now dominates $\mathcal V\times\mathcal V$ by evaluations at $0,\infty$, after reparametrization. At a general map, generic smoothness implies that sections of the pullback tangent bundle span the two fibers jointly. This implies that, in its splitting on $\mathbb{P}^1$,  every summand has positive degree. We therefore have a very free rationalcurve, which we use as the ``handle'' in the construction below.

We now proceed as in Section~\ref{sec:completion}. The klt log Fano perturbation $(X,\Delta+T)$ constructed above and $\mathbb Q$-factoriality implies finite generation of the Cox ring  \cite[Corollary~1.3.2]{BCHM}. Thus, there are nonzero effective Cartier divisors $A_1,\ldots,A_{\mathsf s}$ that generate the pseudoeffective cone after rescaling their classes. For each $j$, apply Proposition \ref{cklt:prop:singular-divisor-test} to $p^*A_j$ to get a family of free curves. Choose one free curve from each family and a very free ``handle'' through a common general point of the trivial-stabilizer locus outside all root divisors. This is possible since all the evaluation images contain dense opens. By Lemma~\ref{lem:smooth-sum}, the resulting curves  through this point have a free smoothing $u$ through the same point whose coarse class is their sum. This construction uses only freeness on the smooth root stack and the common point off its boundary divisor.

Attach $u$ there to the very free handle and place two distinct marks $p_0,q_0$ elsewhere on the handle. For the resulting map $F:C\to\mathcal V$ (the ``comb''), the normalization sequence gives
\[
 H^1(C,F^*T_{\mathcal V}(-p_0-q_0))=0.
\]
Indeed, the handle has vanishing $H^1$ after the twist by $-p_0-q_0$, and the free tooth has vanishing $H^1$ and spans the tangent space at the node. Thus, the pointed smoothing argument in the proof of Lemma~\ref{lem:smooth-sum} smooths the node. Semicontinuity of the above vanishing of $H^1$ impies that we have a very free curve. Its projected class in $X$ is the sum of the original classes, each nonnegative against every $A_j$ by dominance, and the $j$-th curve class has positive $A_j$-degree. Since this class is positive on the finitely many generators of the pseudoeffective one, it is in the interiror of $\operatorname{Nef}_1(X)$.
\end{proof}

\section{Constructions related to divisorial contractions}\label{cklt:sec:divisorialcontraction}
In this section, we work out some technical constructions needed in order to descend a very free rational curve to from a log resolution to $X$. 

\subsection{Cancelling contacts with an exceptional divisor}\label{cklt:sec:cancellation}

To descend to $X$, we must deform a curve off the exceptional divisors of the log resolution while keeping control of two evaluations. We do this by attaching one tooth at each transverse contact. The following proposition explains the choice of its normal degree.

\begin{proposition}[Cancelling transverse contacts]\label{cklt:prop:exceptional-cancellation}
Let $Z$ be a smooth connected separated Deligne--Mumford stack of finite type over $\mathbb{C}$ with quasi-projective coarse space, and let $E_1,\ldots,E_{\mathsf s}$ be disjoint smooth integral Cartier divisors. Suppose each $E_j$ is covered by free maps $u_j:\mathbb{P}^1\to E_j$ of normal degree $-d_j<0$. Fix positive multiples $r_j$ of $d_j$. If the stack $S$ obtained by taking the $r_j$-th root along each $E_j$ admits a very free map from $\mathbb{P}^1$, then $Z\setminus\bigcup_jE_j$ admits a family of maps from $\mathbb{P}^1$ with dominant two-point evaluation.
\end{proposition}

\begin{proof}
Write $G_j$ for the reduced root divisor over $E_j$. On the free locus, the evaluation from the universal curve is smooth. The inverse image of $G_j$ under this evaluation map is smooth, and generic smoothness of the projection to the parameter space gives transversality for a general very free ``handle'' $h:C_0=\mathbb{P}^1\to S$. The inverse image of a codimension-$2$ closed subset under evaluation has codimension at least $2$. Since the universal curve has relative dimension $2$, its projection to the parameter space is not dense. We may therefore take every contact $x_i$ over the open subset of $E_j$ covered by the ``teeth'' $u_j:\PP^1\to E_j$. Choose two distinct marks $p,q$ on the handle off the boundary. If there are no contacts, very freeness gives smooth two-point evaluation, and properness of $C_0$ makes avoidance of the boundary an open condition.

Suppose there are $k>0$ contacts. At $x_i$ over $E_j$, precompose a free tooth through its image with a cover of degree $r_j/d_j$. The resulting map $v_i:C_i=\mathbb{P}^1\to E_j$ has normal line of degree $-r_j$. Since $$G_j\to E_j$$ 
is the gerbe of $r_j$-th roots of $N_{E_j/Z}$, a choice of root $\mathcal O_{C_i}(-1)$ gives a lift $$t_i:C_i\to G_j$$ 
with $t_i^*N_{G_j/S}\simeq\mathcal O_{C_i}(-1)$. Choose $y_i\in C_i$ over the attachment value and an isomorphism $t_i(y_i)\simeq h(x_i)$ in the root gerbe. Identify $y_i$ with $x_i$ to form the comb. It is the Ferrand pushout of its normalization and its node set along the preimages of the nodes. The $2$-pushout property thus glues the component maps to $$F:C\to S,$$ 
see e.g. \cite[Theorem~4.2 and Lemma~4.5]{cklt:AHHR}.

We will smooth the nodes while prescribing the two evaluations. Put $\mathfrak V=F^*T_S(-p-q)$ and $\mathfrak A=\mathfrak V|_{C_0}$. On each tooth, we have
\[
 0\longrightarrow v_i^*T_{E_j}\longrightarrow t_i^*T_S
        \longrightarrow\mathcal O_{C_i}(-1)\longrightarrow0.
\]
Here the gerbe is \'etale locally a product with $B\mu_{r_j}$. The left term is globally generated. Both outer terms have vanishing $H^1$, and the normal line has no global sections. Thus $$H^1(C_i,t_i^*T_S)=0,$$ 
and values of sections at the node span exactly $T_{G_j,h(x_i)}$. 
Very freeness of the handle gives $H^1(C_0,\mathfrak A)=0$. Taking cohomology of the normalization sequence therefore gives, with $N_i=N_{G_j/S,h(x_i)}$,
\begin{equation}\label{cklt:eq:cancellation-normal-cokernel}
 H^1(C,\mathfrak V)=\operatorname{coker}\!\left(H^0(C_0,\mathfrak A)\longrightarrow\bigoplus_i N_i\right).
\end{equation}
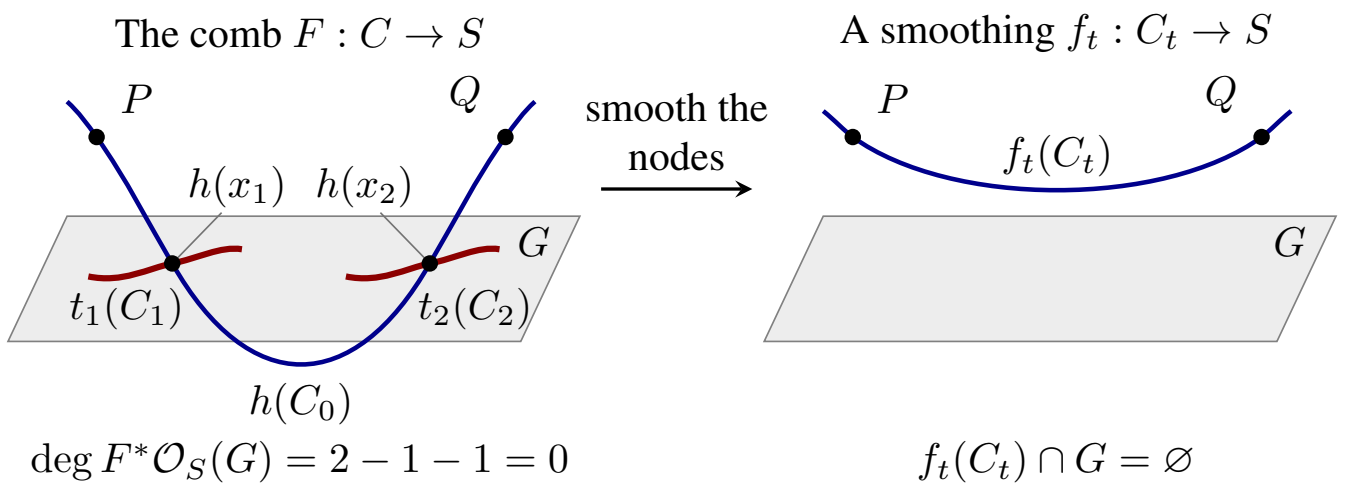
\begin{figure}
    \centering
    \resizebox{0.50\textwidth}{!}{%
\begin{tikzpicture}[x=1cm,y=1cm,font=\footnotesize,>=stealth,
  curve/.style={draw=blue!55!black,line width=1.1pt},
  tooth/.style={draw=red!55!black,line width=1.5pt},
  sheet/.style={draw=black!50,fill=black!7,thin},
  dot/.style={circle,fill=black,inner sep=1.4pt}]

  \foreach \s in {0,6.4}{
    \begin{scope}[xshift=\s cm]
      \draw[sheet] (-0.05,-.08)--(4.29,-.08)--(4.79,.98)--(.45,.98)--cycle;
      \node at (4.40,.76) {$G$};
    \end{scope}
  }

  \node at (2.42,2.53) {The comb $F:C\to S$};
  \coordinate (A) at (1.34,.58);
  \coordinate (B) at (3.52,.58);
  \draw[curve]
    (.45,1.95) .. controls (.55,1.86) and (.62,1.76) .. (.70,1.65)
    .. controls (.95,1.30) and (1.15,.90) .. (A)
    .. controls (1.99,-.56) and (2.87,-.56) .. (B)
    .. controls (3.71,.90) and (3.91,1.30) .. (4.16,1.65)
    .. controls (4.24,1.76) and (4.31,1.86) .. (4.41,1.95);
  \node[below] at (2.43,-.28) {$h(C_0)$};

  \foreach \a/\i/\dx in {1.34/1/.55,3.52/2/-.55}{
    \begin{scope}[xshift=\a cm]
      \draw[thin,black!55] (0,.58)--({.76*\dx},1.01);
      \draw[tooth]
        (-.71,.47) .. controls (-.42,.41) and (-.20,.54) .. (0,.58)
        .. controls (.21,.62) and (.41,.73) .. (.59,.70);
      \node[dot] at (0,.58) {};
      \node at (\dx,1.23) {$h(x_{\i})$};
      \node at ({-.70*\dx},.20) {$t_{\i}(C_{\i})$};
    \end{scope}
  }

  \draw[->,thick] (4.97,1.21)--(6.25,1.21)
    node[midway,above,align=center] {smooth the\\nodes};

  \begin{scope}[xshift=6.4cm]
    \node at (2.42,2.53) {A smoothing $f_t:C_t\to S$};
    \draw[curve]
      (.45,1.87) .. controls (.55,1.79) and (.62,1.71) .. (.70,1.65)
      .. controls (1.53,1.05) and (3.33,1.05) .. (4.16,1.65)
      .. controls (4.24,1.71) and (4.31,1.79) .. (4.41,1.87);
    \node at (2.43,1.49) {$f_t(C_t)$};
  \end{scope}

  \foreach \s in {0,6.4}{
    \begin{scope}[xshift=\s cm]
      \node[dot,label=above right:$P$] at (.70,1.65) {};
      \node[dot,label=above left:$Q$] at (4.16,1.65) {};
    \end{scope}
  }

  \node at (2.42,-1.10) {$\deg F^*\mathcal O_S(G)=2-1-1=0$};
  \node at (8.82,-1.10) {$f_t(C_t)\cap G=\varnothing$};
\end{tikzpicture}%
}
\caption[Cancelling transverse contacts with a root divisor]{\small Cancelling two transverse contacts with the reduced root divisor $G\subset S$. On the left, the handle $h(C_0)$ meets $G$ at
$h(x_1)$ and $h(x_2)$. At each contact, we attach a tooth $t_i:C_i\to G$ with $t_i^*N_{G/S}\simeq\mathcal O_{C_i}(-1)$. Thus the handle contributes $2$ and the two teeth contribute $-1$ each to $\deg F^*\mathcal O_S(G)$. On the right, smoothing both nodes gives a rational curve through the same marked values $P$ and $Q$. Its intersection number with $G$ is $0$, and it contains $P\notin G$, so it is disjoint from $G$.}
\label{fig:cancelling-transverse-contacts}
\end{figure}

Next, we vary the attachment points to give the normal values appearing in this cokernel.

Put two auxiliary markings on each tooth. In an \'etale chart $U$ of the relevant moduli of stable pointed curves, let $u\in U$ be the resulting comb. The node parameters give a smooth map $$\tau:U\to\mathbb{A}^k.$$ 
Its fiber at $0\in \mathbb{A}^k$ parametrizes deformations with the nodes unsmoothed. The $3$-pointed teeth are rigid. On the handle, $T_0=T_{C_0}(-p-q)\simeq\mathcal O_{C_0}$, so the exact sequence for evaluation at the nodes gives
\begin{equation}\label{cklt:eq:cancellation-handle-cokernel}
 T_{U/\mathbb{A}^k,u}=H^1(C_0,T_0(-\textstyle\sum_i x_i))
       =\operatorname{coker}\!\left(H^0(C_0,T_0)\longrightarrow\bigoplus_i T_{x_i}C_0\right).
\end{equation}
To describe the map between the two cokernels, we keep the component maps fixed and move $x_i$ with velocity $w_i$. The two branches now differ at the node by $dh_{x_i}(w_i)$. Transversality makes the resulting linear map $T_{x_i}C_0\to N_i$ an isomorphism for each $i$. A global vector field in $H^0(C_0,T_0)$ gives a section of $\mathfrak A$ by $dh$. Thus the node matchings in the normalization sequence induce a surjection
\begin{equation}\label{cklt:eq:cancellation-equisingular}
 \delta:T_{U/\mathbb{A}^k,u}\twoheadrightarrow H^1(C,\mathfrak V).
\end{equation}

Let $C_U\to U$ be the universal curve and put $\mathcal M_U=\underline{\mathrm{Hom}}_U(C_U,S\times U)$. $C_U$ is proper, flat and finitely presented over $U$, and $S\times U$ is separated of finite type with finite stabilizers. Thus $\mathcal M_U$ is algebraic and locally of finite presentation, see e.g.  \cite[Theorem~1.2(i)--(ii)]{cklt:HR}. We claim that
\[
 (\tau,\operatorname{ev}_p,\operatorname{ev}_q):\mathcal M_U\longrightarrow\mathbb{A}^k\times S\times S
\]
is smooth at $F$. Consider a deformation of $F$ over the quotient of a small Artin extension with kernel $I$, and prescribe lifts of the smoothing parameters and the evaluations (i.e. images in $\mathbb{A}^k\times S\times S$). By smoothness of $\tau$, we have a lift $u'$ of $u\in U$. For the lift $u'$, the gluing obstruction for maps is in $H^1(C,\mathfrak V)\otimes I$, since local lifts with the prescribed values differ by sections vanishing at $p,q$, see e.g. \cite[Theorem~1.5]{Olsson}. Varying $u'$ by $\xi\in T_{U/\mathbb{A}^k,u}\otimes I$ leaves the chosen smoothing parameters unchanged. The obstruction varies linearly with the square-zero extension class, see e.g. \cite[Theorem~6.1 and Remark~6.2]{cklt:Kuhn}. Since the maximal ideal annihilates $I$, the variation is computed on the closed fiber. The change in matching gives $(\delta\otimes 1_I)(\xi)$, so \eqref{cklt:eq:cancellation-equisingular} allows us to cancel the obstruction. This proves the lifting criterion for small extensions, and hence smoothness at $F$ by local finite presentation.

The image of a smooth neighborhood of $F$ contains an open neighborhood of $(0,h(p),h(q))\in \mathbb{A}^k\times S\times S$. Intersect it with $\{\tau_1\cdots\tau_k\ne0\}\subset \mathbb{A}^k\times S\times S$ and the locus where the value at $p$ misses the boundary divisor. Its projection to $S\times S$ contains a nonempty open, so two-point evaluation for the resulting family is dominant. For each $G_j$, the degree of $F^*\mathcal O_S(G_j)$ is $0$, since every contact contributes $1$ on the handle and $-1$ on its tooth. The degree is constant in the smoothing. On each smooth fiber the divisor section is nonzero at $p$ and hence nowhere zero. The maps therefore lie in the complement of $G_j$. The two evaluation markings and one auxiliary marking trivialize the smooth rational source family, giving maps from a fixed $\mathbb{P}^1$. A general pair of evaluation values has distinct coarse images.
\end{proof}


\subsection{Adjunction along the exceptional divisor}\label{cklt:sec:adjunction}
At the contraction of an exceptional divisor, we need curves in the exceptional divisor to serve as teeth in Proposition~\ref{cklt:prop:exceptional-cancellation}. Adjunction determines the root orders on this divisor, and induction on dimension gives the teeth.

For a plt pair $(Y,E+B)$ with $E$ prime of coefficient $1$, $E$ is normal. Taking residues of sufficiently divisible log pluricanonical forms along $E$ gives the different $\operatorname{Diff}_E B$ and the adjunction formula
\[
 (K_Y+E+B)|_E\sim_{\mathbb Q} K_E+\operatorname{Diff}_E B .
\]

\begin{lemma}[Orbifold adjunction]\label{cklt:lem:orbifold-adjunction}
Let $Y$ be normal and $\mathbb Q$-factorial, and let $(Y,E+B)$ be plt with $E$ prime and $B=\sum_i(1-1/m_i)B_i$ of finite standard coefficients, where $E$ is distinct from every $B_i$. Then there is an open subvaritey $$Y^{\mathrm{orb}}\subset Y$$ 
containing every point of codimension at most $2$, and a smooth separated Deligne--Mumford stack $\mathcal P$ with coarse space $Y^{\mathrm{orb}}$, trivial generic stabilizer and divisorial branch coefficients given by those of $B$. The reduced divisor $\mathcal E\subset \mathcal{P}$ over $E$ is smooth Cartier, has trivial generic stabilizer and coarse space $E\cap Y^{\mathrm{orb}}$. Its divisorial branch coefficients on $E$ are those of $\operatorname{Diff}_E B$. In particular, the different has standard coefficients.
\end{lemma}

\begin{proof}
Choose $a_i>0$ such that $a_iB_i$ Cartier. Consider the fiber product over $Y$ of the root stacks $\sqrt[a_i m_i]{(Y,a_iB_i)}$ and its normalization; if $B=0$, this is just $Y$. Over a trivializing open $\mathsf U$, a chart is the normalization $\mathsf U'$ of the finite cover adjoining these roots. Its invariant algebra under the product of the groups of roots of unity is finite over $\mathcal O_{\mathsf U}$, and agrees with $\mathcal O_{\mathsf{U}}$ away from  the boundary divisor. Since the spaces are normality, these algebras are equal. Over the geometric generic point of $B_i$, put $a=a_i$, $m=m_i$, and let $x$ be a uniformizer. After extracting a root of a unit \'etale locally, the equation is $z^{am}=x^a$. Its normalized branches $z^m=\zeta x$, $\zeta\in\mu_a$, are permuted transitively by $\mu_{am}$ with stabilizer $\mu_m$. Thus the ramification index is $m$. The other covers are \'etale there.

The resulting finite normal charts $g:\mathsf U'\to \mathsf U$ satisfy $K_{\mathsf U'}=g^*(K_{\mathsf U}+B|_{\mathsf U})$ by the ramification formula. Since $(\mathsf U,B|_{\mathsf U})$ is klt, $\mathsf U'$ is also klt, see  \cite[Proposition~5.20]{cklt:KM}. Therefore, by \cite[Lemma~14]{cklt:CGG}, $\mathsf{U}'$ has quotient singularities away from a subset of codimension at least $3$. \'Etale locally, we can present $\mathsf{U}'$ as quotients by finite groups. By linearization of the finite group action on the completed smooth local ring and Artin approximation \cite[Corollary~2.6]{cklt:Artin}, these quotient charts are algebraic. Remove the images of the nonquotient loci and take the canonical smooth stacks of the remaining charts. By functoriality under \'etale restriction and the local group actions, this gives a smooth separated stack $\mathcal P$ over the resulting open of $Y$. Its generic stabilizer is trivial and the indices in codimension $1$ are unchanged. On the basic open for $B$, in the sense of (\ref{eqn:basic_open}), the characterization in \cite[Theorem~1]{cklt:GS} identifies it with the given product of root stacks.

Take, \'etale locally, a smooth finite uniformizing chart $v:\widehat{\mathsf{U}}\to \mathsf{U}$ for $\mathcal P$. It is unramified over the generic point of $E$, so $\widehat E=v^*E$ is reduced and
\[
 K_{\widehat{\mathsf U}}+\widehat E=v^*(K_{\mathsf U}+E+B).
\]
It follows from comparing discrepancies that $(\widehat{\mathsf U},\widehat E)$ is plt. Its coefficient-$1$ boundary divisor components are disjoint and normal. Removing the images of their singular loci affects only codimension at least $3$ in $Y$ and leaves a smooth Cartier divisor $\mathcal E$. Since taking invariants under finite group actions is exact, the quotient of $\widehat E$ can be identified with its reduced image in $\mathsf U$, giving the coarse space as claimed. Its generic stabilizer is trivial. Taking residues in the above equality, we get 
$$K_{\widehat E}=v_E^*(K_E+\operatorname{Diff}_E B),$$ 
where $v_E$ is the restriction of $v$. If the index over a prime of $E$ is $e$, the ramification coefficient upstairs is $e-1$, so the coefficient of the different is $1-1/e$. Every codimension-$1$ point of $E$ remains in this calculation.
\end{proof}

\subsection{Curves in a contracting fiber}\label{cklt:subsec:vertical-family}

Let $$\phi:Y\to Y_1$$ 
be a projective divisorial extremal contraction with exceptional divisor $E$ and center $S$. Assume that we are in the situation of Lemma~\ref{cklt:lem:orbifold-adjunction}, and suppose that both $-(K_Y+E+B)$ and $-E$ are relatively ample. Suppose Theorem~\ref{cklt:thm:dlt-descent} holds in lower dimensions. We now construct a covering family of free curves on $\mathcal E$, contracted by $\phi$ and of fixed negative normal degree in $\mathcal P$.

Write $$E\to S'\to S$$ 
for the Stein factorization. The generic fiber over $S'$ is normal because $E$ is normal, and geometrically integral because $\mathbb{C}(S')$ is algebraically closed in $\mathbb{C}(E)$ in characteristic zero. It is geometrically normal since the field $\mathbb{C}(S')$ is perfect. By generic flatness and openness of these properties, we have a normal integral general fiber $F$ over $S'_{\mathrm{reg}}$, with $0<\dim F<\dim Y$. Set $B_F=(\operatorname{Diff}_E B)|_F$. By adjunction, $(E,\operatorname{Diff}_E B)$ klt. Take a log resolution of this pair and restrict it to $F$. By generic smoothness on the resolution and its strata, we have a log resolution of $F$ with the same log pullback coefficients on components horizontal with respect to $E\to S'$. Vertical components miss $F$. Distinct horizontal primes have reduced restrictions at their generic points and have no common divisorial component on $F$. Thus $(F,B_F)$ is klt with finite standard coefficients, and adjunction gives $$K_F+B_F\sim_{\mathbb Q}(K_Y+E+B)|_F,$$ 
whose negative is ample.

On the intersection of $E\cap Y^{\mathrm{orb}}$ with the basic open for the different, $\mathcal E$ is its given root stack by \cite[Theorem~1]{cklt:GS}. The complement has codimension at least $2$ in $E$ and hence, by the fiber dimension theorem, also on a general fiber which it meets. After shrinking $S'$, the open and its disjoint smooth boundary components are smooth over $S'$. Since formations of root stacks and tame coarse spaces commute with base change (see e.g. \cite{cklt:AOV}), on this common open in $F$, the restriction $\mathcal E_F$ agrees with the root stack for $B_F$.

By induction on dimension, we get a very free map on the basic-open root stack of $(F,B_F)$. By smoothness of universal evaluation on the free locus, a general deformation avoids the omitted codimension-$2$ subset. Denote the resulting map to $\mathcal E_F$ by $v:\PP^1\to \mathcal{E}_F$. Along $v$, the relative tangent sequence has ample subbundle $v^*T_{\mathcal E_F}$ and quotient $\mathcal O_{\mathbb{P}^1}^{\oplus\dim S'}$. The sequence splits because  $H^1(\PP^1, v^*T_{\mathcal E_F})=0$. Hence $v$ is free in $\mathcal E$ and its deformations give a covering family. Its deformations remain vertical with respect to $E\to S'$ and has degree $0$ against the pullback of an ample divisor on $Y_1$. The divisor class of $\mathcal E$ is the rational pullback of $E$ since the coarse map is unramified at its generic point. By the assumption of relative ampleness of $-E$, the negative normal degree of $v$ is negative. This degree, denoted by $-d$, is an integer and constant in the family.

\section{Descending a very free curve}\label{cklt:sec:descent}
Here we work out the process of descending a rational curve from a log resolution to $X$.

\begin{proof}[Proof of Theorem~\ref{cklt:thm:dlt-descent}]
We proceed by induction on $\dim X$. For dimension $1$, $X$ is smooth and Proposition~\ref{cklt:prop:initial} applied to the identity resolution gives the desired very free map of interior class. We verify the contact assertions below in all dimensions.

Suppose $\text{dim}\, X>1$. We take a small projective $\mathbb Q$-factorialization $b:X'\to X$ \cite[Corollary~1.4.3]{BCHM}. For the strict transform $\Delta'$ we have 
$$K_{X'}+\Delta'=b^*(K_X+\Delta),$$
so $(X',\Delta')$ is klt and weak log Fano. The basic-open root stacks agree away from codimension at least $2$ on both sides. A general very free map on $X'$ avoids the omitted locus and induces one to $X$. Its class stays interior since $b^*$ embeds $N^1(X)_{\mathbb R}$ into $N^1(X')_{\mathbb R}$ and preserves pseudoeffectiveness of divisor classes. By projection formula, the pairing of this curve class with every nonzero pseudoeffective class on $X$ is positive. We may thus work on the $\mathbb Q$-factorial model $(X',\Delta')$.

\smallskip\noindent\emph{The relative MMP.}
Take a projective log resolution $$p:W\to X.$$ 
Write $D_W$ for the strict transform of $\Delta$ and write $F_W=\sum_jF_j$ for the reduced exceptional divisor. Since $(X,\Delta)$ is klt, we have
\[
 K_W+D_W+F_W=p^*(K_X+\Delta)+\sum_j(1+a(F_j,X,\Delta))F_j
\]
which has effective exceptional part with support $F_W$. The relative dlt minimal model program with ample scaling gives a finite sequence of birational models that contracts all $F_j$, see \cite[Theorem~3.5]{cklt:Birkar}, through $\mathbb Q$-factorial dlt models $(Y,D_Y+F_Y)$. The final step gives a small morphism to $X$, which is an isomorphism. Indeed, a relatively ample Cartier divisor would be equal to the rational pullback of its pushforward by $\mathbb Q$-factoriality of $X$ and smallness, hence has degree $0$ on any contracted curve.

Fix this sequence of MMP. By strict negativity on its finitely many divisorial rays and continuity in the coefficients, there is $R\ge2$ such that all these degrees remain negative when the coefficients of $F_Y$ vary independently in $[1-1/R,1]$. On the basic open for $D_Y+F_Y$, write $\mathcal Y^{\circ}(\boldsymbol\ell)$ for the stack with the given roots along $D_Y$ and orders $\ell_j$ along the remaining exceptional primes. For each model $Y$ in this sequence, we prove that
\[{\bf (\star)}\quad
 \begin{gathered}
 \text{for every tuple }\boldsymbol\ell\text{ with all }\ell_j\ge R,
 \text{ there is a very free map }\PP^1\to\mathcal Y^{\circ}(\boldsymbol\ell)\\
 \text{ whose coarse class projects to the interior of } \operatorname{Nef}_1(X).
 \end{gathered}
\]
For $W$, we apply Proposition~\ref{cklt:prop:initial} and deform away from intersections of boundary components. When a step goes through a flip, the basic-open root stacks agree on the common isomorphism locus up to codimension $2$. A general deformation there preserves the tangent bundle and the class over $X$.

\smallskip
\noindent\emph{A divisorial contraction.}
Suppose $\phi:Y\to Y_1$ contracts $E$ onto $S$. Fix the tuple on $Y_1$ and add the strict transforms of the surviving exceptional primes with coefficients $1-1/\ell_j$ to $D_Y$; this gives $B$. The pair $(Y,E+B)$ is plt, since lowering coefficients preserves dlt and its only coefficient-$1$ stratum is $E$. By the choice of $R$, $-(K_Y+E+B)$ is $\phi$-ample. By \cite[Lemma 3.39]{cklt:KM}, $E$ has negative degree on the contracted ray, so $-E$ is also $\phi$-ample. Every contracted curve lies in $E$, and $Y\setminus E\simeq Y_1\setminus S$.

By Lemma~\ref{cklt:lem:orbifold-adjunction} and Section~\ref{cklt:subsec:vertical-family}, we have $\mathcal E\subset\mathcal P$ and a covering vertical family of normal degree $-d$. Now, we choose $r\ge R$ divisible by $d$. Apply ${\bf (\star)}$ on $Y$ with root order $r$ on $E$ and the fixed orders on the other components. On the common basic open, its root stack agrees with $\sqrt[r]{(\mathcal P,\mathcal E)}$. A general deformation avoids the omitted codimension-$2$ locus and gives a very free handle there.

By Proposition~\ref{cklt:prop:exceptional-cancellation}, we can smooth the comb, made from the handle and the vertical teeth, to maps in $\mathcal P\setminus\mathcal E$ with dominant two-point evaluation. By the argument in Proposition~\ref{cklt:prop:initial}, the general map is very free. Its general deformation avoids the codimension-$2$ complement of the common basic open in $Y\setminus E\simeq Y_1\setminus S$, where the roots along the surviving components coincide. Consequently, this very free rational  curve induces one on $\mathcal Y_1^{\circ}(\boldsymbol\ell)$. Every tooth is contracted over $X$, and degrees of pullbacks from $X$ are constant in the smoothing. Thus the projected class is that of the handle. This proves ${\bf (\star)}$ for the arbitrary tuple on $Y_1$ and finishes the induction along the MMP.

\smallskip
\noindent\emph{The coarse map and its contacts.}
Let $$g:C=\mathbb{P}^1\to\mathcal X$$ 
be a general deformation of the very free curve we just constructed, and put $f=\pi\mathbin{\!\scriptstyle\circ\!} g$. On the small $\mathbb Q$-factorialization, the strict transform of each $D_i$ has nonzero effective numerical class, hence its degree on the chosen interior curve class is positive. The general curve meets the complement of the boundary divisor. It therefore meets each strict transform, and transfers with these intersections to the common basic open of $X$.

We apply the contact calculation in Lemma~\ref{lem:free-contacts} and its proof on $X^{\circ}$. It applies on this smooth open because its ingredients are smooth evaluation of free maps, properness of the source $C$, and the root charts $x=z^{m_i}$ along the boundary. It gives transversality upstairs, contact order $m_i$ downstairs, and a family of marked coarse maps with fixed contact data after shrinking and labeling the contacts by finite \'etale base change. By very freeness, $H^1(C,g^*T_{\mathcal X}(-2))=0$, so the evaluation at two distinct points is smooth. Dominance continues to hold on the coarse space and under the indicated base change.

For logarithmic very freeness, we use the comparison $N_f^{\log}\simeq [T_C\to g^*T_{\mathcal X}]$ in the proof of Lemma~\ref{lem:free-contacts}, with the two terms in degrees $-1,0$. At each transverse contact, the quotient of the ordinary complex by the logarithmic one is the isomorphism from the source tangent line to the root normal line. Twisting the comparison by $\mathcal O_C(-2)$ gives $\mathbb H^1(C,N_f^{\log}(-2))=0$ by the preceding vanishing and $H^2(C,T_C(-2))=0$. This is logarithmic very freeness in the sense of \cite[Definition~4.2]{CLT}. The comparison also gives $\deg N_f^{\log}=-(K_X+\Delta)\cdot f_*[C]-2$, since $K_{\mathcal X}=\pi^*((K_X+\Delta)|_{X^{\circ}})$.
\end{proof}

\end{document}